\documentclass[11pt]{article}%
\usepackage{amsmath}
\usepackage{amsfonts}
\usepackage{amssymb}
\usepackage{graphicx}
\usepackage{a4wide}%
\providecommand{\U}[1]{\protect\rule{.1in}{.1in}}
\newtheorem{theorem}{Theorem}

\newtheorem{corollary}[theorem]{Corollary}

\newtheorem{example}[theorem]{Example}

\newtheorem{lemma}[theorem]{Lemma}

\newtheorem{proposition}[theorem]{Proposition}
\newtheorem{remark}[theorem]{Remark}

\begin{document}

\title{On quadratic optimization problems and canonical duality theory}
\author{C. Z\u{a}linescu\\Faculty of Mathematics, University Al. I. Cuza Iasi, Iasi, Romania}
\maketitle

Canonical duality theory (CDT) is advertised by its author DY Gao as
\textquotedblleft a breakthrough methodological theory that can be used not
only for modeling complex systems within a unified framework, but also for
solving a large class of challenging problems in multidisciplinary fields of
engineering, mathematics, and sciences."

DY Gao solely or together with some of his collaborators applied CDT for
solving some quadratic optimization problems with quadratic constraints.
Unfortunately, in almost all papers we read on CDT there are unclear
definitions, non convincing arguments in the proofs, and even false results.

The aim of this paper is to treat rigorously quadratic optimization problems
by the method suggested by CDT and to compare what we get with the results
obtained by DY Gao and his collaborators on this topic in several papers.

\section{Notations and preliminary results\label{sec1}}

Let us consider the quadratic functions $q_{k}:\mathbb{R}^{n}\rightarrow
\mathbb{R}$ for $k\in\overline{0,m}$, that is $q_{k}(x):=\tfrac{1}%
{2}\left\langle x,A_{k}x\right\rangle -\left\langle b_{k},x\right\rangle
+c_{k}$ for $x\in\mathbb{R}^{n}$ with given $A_{k}\in\mathfrak{S}_{n}$,
$b_{k}\in\mathbb{R}^{n}$ (seen as column vector) and $c_{k}\in\mathbb{R}$ for
$k\in\overline{0,m}$, where $\mathfrak{S}_{n}$ denotes the class of symmetric
matrices from $\mathfrak{M}_{n}:=\mathbb{R}^{n\times n}$, and $\left\langle
\cdot,\cdot\right\rangle $ denotes the usual inner product on $\mathbb{R}^{n}%
$. For $k\in\mathbb{N}^{\ast}$ We set%
\[
\mathbb{R}_{+}^{k}:=\{\eta\in\mathbb{R}^{k}\mid\eta_{i}\geq0~\forall
i\in\overline{1,k}\},\quad\mathbb{R}_{-}^{k}:=-\mathbb{R}_{+}^{k}%
,\quad\mathbb{R}_{++}^{k}:=\operatorname*{int}\mathbb{R}_{+}^{k}%
,\quad\mathbb{R}_{--}^{k}:=-\mathbb{R}_{++}^{k}.
\]

The fact that $A\in\mathfrak{S}_{n}$ is positive (semi) definite is denoted by
$A\succ0$ $(A\succeq0)$ and we set $\mathfrak{S}_{n}^{+}:=\{A\in
\mathfrak{S}_{n}\mid A\succeq0\}$, $\mathfrak{S}_{n}^{++}:=\{A\in
\mathfrak{S}_{n}\mid A\succ0\}$ ; it is well known that $\mathfrak{S}_{n}%
^{++}=\operatorname*{int}\mathfrak{S}_{n}^{+}$. In this paper we consider
quadratic minimization problems with (quadratic) equality and inequality
constraints. With this aim, we fix a set $J\subset\overline{1,m}$
corresponding to the equality constraints; the set $J^{c}:=\overline
{1,m}\setminus J$ will correspond to the inequality constraints. So, the
general problem is

\medskip

$(P_{J})$ $~~\min$ $q_{0}(x)$ ~s.t. $x\in X_{J}$, \medskip

\noindent where
\[
X_{J}:=\{x\in\mathbb{R}^{n}\mid\left[  \forall j\in J:q_{j}(x)=0\right]
~\wedge~\left[  \forall j\in J^{c}:q_{j}(x)\leq0\right]  \}.
\]
For later use we introduce also the set
\[
\Gamma_{J}:=\{(\lambda_{1},...,\lambda_{m})\in\mathbb{R}^{m}\mid\lambda
_{j}\geq0~\forall j\in J^{c}\}. \label{r-gamma}%
\]
Clearly, for $J=\overline{1,m}$ $(P_{J})$ becomes the quadratic minimization
problem with (quadratic) equality constraints denoted $(P_{e})$ with
$X_{e}:=X_{\overline{1,m}}$ its feasible set, while for $J=\emptyset$
$(P_{J})$ becomes the quadratic minimization problem with inequality
constraints denoted $(P_{i})$ with $X_{i}:=X_{\emptyset}$ its feasible set.
Clearly $X_{e}\subset X_{J}\subset X_{i}$, the inclusions being strict in
general when $\emptyset\neq J\neq\overline{1,m}$. Observe that any
optimization problem with equality constraints can be seen as a problem with
inequality constraints because the equality constraint $h(x)=0$ can be
replaced by the inequality constraints $g_{1}(x):=h(x)\leq0$ and
$g_{2}(x):=-h(x)\leq0$. Excepting linear programming, such a procedure is not
used in general because the constraints qualification conditions are very
different for problems with equality constraints and those with inequality constraints.

To the family $(q_{k})_{k\in\overline{0,m}}$ we associate the Lagrangian
$L:\mathbb{R}^{n}\times\mathbb{R}^{m}\rightarrow\mathbb{R}$ defined by
\[
L(x,\lambda):=q_{0}(x)+\sum\nolimits_{j=1}^{m}\lambda_{j}q_{j}(x)=\tfrac{1}%
{2}\left\langle x,A(\lambda)x\right\rangle -\left\langle x,b(\lambda
)\right\rangle +c(\lambda), \label{r-L}%
\]
where $A(\lambda)x:=[A(\lambda)]\cdot x$ and%
\[
A(\lambda):=\sum\nolimits_{k=0}^{m}\lambda_{k}A_{k},\quad b(\lambda
):=\sum\nolimits_{k=0}^{m}\lambda_{k}b_{k},\quad c(\lambda):=\sum
\nolimits_{k=0}^{m}\lambda_{k}c_{k},
\]
with $\lambda_{0}:=1$ and $\lambda:=(\lambda_{1},...,\lambda_{m})^{T}%
\in\mathbb{R}^{m}$. Clearly, $A:\mathbb{R}^{m}\rightarrow\mathfrak{S}_{n}$,
$b:\mathbb{R}^{m}\rightarrow\mathbb{R}^{n}$, $c:\mathbb{R}^{m}\rightarrow
\mathbb{R}$ defined by the above formulas are affine mappings.

Moreover, one considers the sets
\begin{gather}
Y_{0}:=\{\lambda\in\mathbb{R}^{m}\mid\det A(\lambda)\neq0\},\label{r-s0}\\
Y^{+}:=\{\lambda\in\mathbb{R}^{m}\mid A(\lambda)\succ0\},\quad Y^{-}%
:=\{\lambda\in\mathbb{R}^{m}\mid A(\lambda)\prec0\}.\label{r-s0pm}%
\end{gather}
Observe that $Y_{0}$ is a (possible empty) open set, while $Y^{+}$ and $Y^{-}$
are (possibly empty) open and convex sets. Sometimes one uses also the sets
\begin{gather}
Y_{\operatorname{col}}:=\{\lambda\in\mathbb{R}^{m}\mid b(\lambda
)\in\operatorname{Im}A(\lambda)\},\label{r-yext}\\
Y_{\operatorname{col}}^{+}:=\{\lambda\in Y_{\operatorname{col}}\mid
A(\lambda)\succeq0\},\quad Y_{\operatorname{col}}^{-}:=\{\lambda\in
Y_{\operatorname{col}}\mid A(\lambda)\preceq0\},\label{r-yextp}%
\end{gather}
where for $F\in\mathbb{R}^{m\times n}$ we set $\operatorname{Im}F:=\{Fx\mid
x\in\mathbb{R}^{n}\}$ and $\ker F:=\{x\in\mathbb{R}^{n}\mid Fx=0\}$. Clearly,
$Y_{0}\subset Y_{\operatorname{col}}$, $Y^{+}\subset Y_{\operatorname{col}%
}^{+}$, $Y^{-}\subset Y_{\operatorname{col}}^{-}$, and $Y_{\operatorname{col}%
}$ is neither open, nor closed (in general). Unlike for $Y^{+}$, the convexity
of $Y_{\operatorname{col}}^{+}$ is less obvious. In fact the next (probably
known) result holds.

\begin{lemma}
\label{lem-im}\emph{(i)} Let $A,B\in\mathfrak{S}_{n}^{+}$. Then
$\operatorname{Im}(A+B)=\operatorname{Im}A+\operatorname{Im}B$.

\emph{(ii)} Let $A\in\mathfrak{S}_{n}$ and $a\in\mathbb{R}^{n}$, and set
$q(x):=\tfrac{1}{2}\left\langle x,Ax\right\rangle -\left\langle
a,x\right\rangle $. Then $q(x_{1})=q(x_{2})$ for all $x_{1},x_{2}\in
\mathbb{R}^{n}$ such that $Ax_{1}=Ax_{2}=a$.
\end{lemma}

Proof. (i) It is known that $\operatorname{Im}F=(\ker F)^{\perp}$, and so
$\mathbb{R}^{n}=\operatorname{Im}F+\ker F$, provided $F\in\mathfrak{S}_{n}$.
Moreover, using Schwarz' inequality for positive semi-definite matrices
(operators) we have that $\ker F=\{x\in\mathbb{R}^{n}\mid\left\langle
x,Fx\right\rangle =0\}$ whenever $F\in\mathfrak{S}_{n}^{+}$. Since
$A+B\in\mathfrak{S}_{n}^{+}$ we get
\begin{align*}
\left(  \operatorname{Im}(A+B)\right)  ^{\perp}  &  =\ker(A+B)=\{x\in
\mathbb{R}^{n}\mid\left\langle x,(A+B)x\right\rangle =0\}\\
&  =\ker A\cap\ker B=\left(  \operatorname{Im}A\right)  ^{\perp}\cap\left(
\operatorname{Im}B\right)  ^{\perp}=\left(  \operatorname{Im}%
A+\operatorname{Im}B)\right)  ^{\perp},
\end{align*}
whence the conclusion.

(ii) Take $x_{1},x_{2}\in\mathbb{R}^{n}$ such that $Ax_{1}=Ax_{2}=a;$ setting
$x:=x_{1}$ and $u:=x_{2}-x_{1}$, we have that $x_{2}=x+u$ and $Au=0$. It
follows that $\left\langle a,u\right\rangle =\left\langle Ax,u\right\rangle
=\left\langle x,Au\right\rangle =0$, and so%
\[
q(x+u)=\tfrac{1}{2}\left\langle x+u,A(x+u)\right\rangle -\left\langle
a,x+u\right\rangle =\tfrac{1}{2}\left\langle x,Ax\right\rangle -\left\langle
a,x\right\rangle =q(x),
\]
whence $q(x_{2})=q(x_{1})$. \hfill$\square$

\begin{corollary}
\label{c-yqi-ycoli}With the previous notations and assumptions,
$Y_{\operatorname{col}}^{+}$ and $Y_{\operatorname{col}}^{-}$ are convex.
Moreover, if $Y^{+}$ (resp.~$Y^{-}$) is nonempty, then $Y^{+}%
=\operatorname*{int}Y_{\operatorname{col}}^{+}$ (resp.~$Y^{-}%
=\operatorname*{int}Y_{\operatorname{col}}^{-}$).
\end{corollary}

Proof. Take $\lambda,\lambda^{\prime}\in Y_{\operatorname{col}}^{+}$ and
$\alpha\in(0,1)$. From the definition of $Y_{\operatorname{col}}^{+}$ and
Lemma \ref{lem-im}~(i), taking into account that $A$ and $b$ are affine, we
get
\begin{align*}
b(\alpha\lambda+(1-\alpha)\lambda^{\prime}) &  =\alpha b(\lambda
)+(1-\alpha)b(\lambda^{\prime})\in\alpha\operatorname{Im}A(\lambda
)+(1-\alpha)\operatorname{Im}A(\lambda^{\prime})\\
&  =\operatorname{Im}[\alpha A(\lambda)]+\operatorname{Im}[(1-\alpha
)A(\lambda^{\prime})]=\operatorname{Im}[\alpha A(\lambda)+(1-\alpha
)A(\lambda^{\prime})]\\
&  =\operatorname{Im}A(\alpha\lambda+(1-\alpha)\lambda^{\prime}),
\end{align*}
and so $\alpha\lambda+(1-\alpha)\lambda^{\prime}\in Y_{\operatorname{col}}%
^{+}$. The proof of the convexity of $Y_{\operatorname{col}}^{-}$ is similar.

Assume now that $Y^{+}\neq\emptyset$ and take $\lambda_{0}\in Y^{+}$,
$\lambda\in Y_{\operatorname{col}}^{+}$ and $\alpha\in(0,1)$. Then
$A(\alpha\lambda_{0}+(1-\alpha)\lambda)=\alpha A(\lambda_{0})+(1-\alpha
)A(\lambda)\succ0$, and so $\alpha\lambda_{0}+(1-\alpha)\lambda\in Y^{+}$.
Taking the limit for $\alpha\rightarrow0$ we obtain that $\lambda
\in\operatorname*{cl}Y^{+}$. Hence $Y^{+}\subset Y_{\operatorname{col}}%
^{+}\subset\operatorname*{cl}Y^{+}$, and so
\[
Y^{+}=\operatorname*{int}Y^{+}\subset\operatorname*{int}Y_{\operatorname{col}%
}^{+}\subset\operatorname*{int}(\operatorname*{cl}Y^{+})=Y^{+}.
\]
The proof is complete. \hfill$\square$

\medskip

Of course, for every $(x,\lambda)\in\mathbb{R}^{n}\times\mathbb{R}^{m}$ we
have that
\begin{equation}
\nabla_{x}L(x,\lambda)=A(\lambda)\cdot x-b(\lambda),\quad\nabla_{xx}%
^{2}L(x,\lambda)=A(\lambda),\quad\nabla_{\lambda}L(x,\lambda)=\left(
q_{j}(x)\right)  _{j\in\overline{1,m}}. \label{r-nxL}%
\end{equation}

Hence $L(\cdot,\lambda)$ is (strictly) convex for $\lambda\in
Y_{\operatorname{col}}^{+}$ $(\lambda\in Y^{+})$ and (strictly) concave for
$\lambda\in Y_{\operatorname{col}}^{-}$ $(\lambda\in Y^{-})$. Moreover, for
$\lambda\in Y_{0}$ we have that $\nabla_{x}L(x,\lambda)=0$ iff $x=[A(\lambda
)]^{-1}\cdot b(\lambda)$, written $A(\lambda)^{-1}b(\lambda)$ in the sequel.

Let us consider now the (dual objective) function%

\begin{equation}
D:Y_{\operatorname{col}}\rightarrow\mathbb{R},\quad D(\lambda):=L(x,\lambda
)\text{ with }A(\lambda)x=b(\lambda);\label{r-pd0}%
\end{equation}
$D$ is well defined because for $x_{1},x_{2}\in\mathbb{R}^{n}$ with
$A(\lambda)x_{1}=A(\lambda)x_{2}=b(\lambda)$, by Lemma \ref{lem-im}~(ii), we
have that $L(x_{2},\lambda)=L(x_{1},\lambda)$. In particular,%
\[
\big[\lambda\in Y_{0}~\text{~}\wedge~~x=\left(  A(\lambda)\right)  ^{-1}\cdot
b(\lambda)\big]\Longrightarrow L(x,\lambda)=D(\lambda).\label{r-lpd}%
\]

Of course
\begin{equation}
D(\lambda)=L\big(A(\lambda)^{-1}b(\lambda),\lambda\big)=-\tfrac{1}%
{2}\big\langle b(\lambda),A(\lambda)^{-1}b(\lambda)\big\rangle+c(\lambda
)\quad\forall\lambda\in Y_{0}.\label{r-pd}%
\end{equation}

\begin{lemma}
\label{lem-qperfdual}Let $(\overline{x},\overline{\lambda})\in\mathbb{R}%
^{n}\times\mathbb{R}^{m}$ be such that $\nabla_{x}L(\overline{x}%
,\overline{\lambda})=0$ and $\left\langle \overline{\lambda},\nabla_{\lambda
}L(\overline{x},\overline{\lambda})\right\rangle =0$. Then $\overline{\lambda
}\in Y_{\operatorname{col}}$ and
\begin{equation}
q_{0}(\overline{x})=L(\overline{x},\overline{\lambda})=D(\overline{\lambda}).
\label{r-qlpd}%
\end{equation}
In particular, $\overline{x}\in X_{e}$ and (\ref{r-qlpd}) hold if
$(\overline{x},\overline{\lambda})$ is a critical point of $L$, that is
$\nabla L(\overline{x},\overline{\lambda})=0$.
\end{lemma}

Proof. Because $0=\nabla_{x}L(\overline{x},\overline{\lambda})=A(\overline
{\lambda})\overline{x}-b(\overline{\lambda})$, it is clear that $\overline
{\lambda}\in Y_{\operatorname{col}}$ and $L(\overline{x},\overline{\lambda
})=D(\overline{\lambda})$ by the definition of $D$. On the other hand,
\[
L(\overline{x},\overline{\lambda})=q_{0}(\overline{x})+\sum\nolimits_{j=1}%
^{m}\overline{\lambda}_{j}q_{j}(\overline{x})=q_{0}(\overline{x})+\left\langle
\overline{\lambda},\nabla_{\lambda}L(\overline{x},\overline{\lambda
})\right\rangle =q_{0}(\overline{x}).
\]
The last assertion follows from the expression of $\nabla_{\lambda}%
L(\overline{x},\overline{\lambda})$ in (\ref{r-nxL}). \hfill$\square$

\medskip

Formula (\ref{r-qlpd}) is related to the so-called \textquotedblleft
complimentary-dual principle\textquotedblright\ (see \cite[p.~NP11]%
{GaoRuaLat:16}, \cite[p.~13]{GaoRuaLat:17}) and sometimes is called the
\textquotedblleft perfect duality formula".

\begin{proposition}
\label{lem-pd}\emph{(i)} The following representation of $D$ holds:
\begin{equation}
D(\lambda)=\left\{
\begin{array}
[c]{ll}%
\min_{x\in\mathbb{R}^{n}}L(x,\lambda) & \text{if }\lambda\in
Y_{\operatorname{col}}^{+},\\
\max_{x\in\mathbb{R}^{n}}L(x,\lambda) & \text{if }\lambda\in
Y_{\operatorname{col}}^{-},
\end{array}
\right.  \label{r-1}%
\end{equation}
the value of $D(\lambda)$ being attained at any $x\in\mathbb{R}^{n}$ such that
$A(\lambda)x=b(\lambda)$ whenever $\lambda\in Y_{\operatorname{col}}^{+}\cup
Y_{\operatorname{col}}^{-};$ in particular, $D(\lambda)$ is attained uniquely
at $x:=A(\lambda)^{-1}b(\lambda)$ for $\lambda\in Y^{+}\cup Y^{-}$.

\emph{(ii)} $D$ is concave and upper semi\-continuous on
$Y_{\operatorname{col}}^{+}$, and convex and lower semi\-continuous on
$Y_{\operatorname{col}}^{-}$.

\emph{(iii)} Let $J\subset\overline{1,m}$ and $(\overline{x},\overline
{\lambda})\in X_{J}\times\mathbb{R}^{m}$ be such that $\nabla_{x}%
L(\overline{x},\overline{\lambda})=0$ and $\left\langle \overline{\lambda
},\nabla_{\lambda}L(\overline{x},\overline{\lambda})\right\rangle =0$. Then
$\overline{\lambda}\in Y_{\operatorname{col}};$ moreover%
\begin{gather*}
\overline{\lambda}\in\Gamma_{J}\cap Y_{\operatorname{col}}^{+}\Longrightarrow
D(\overline{\lambda})=\max\left\{  D(\lambda)\mid\lambda\in\Gamma_{J}\cap
Y_{\operatorname{col}}^{+}\right\}  ,\label{r-2a}\\
\overline{\lambda}\in(-\Gamma_{J})\cap Y_{\operatorname{col}}^{-}%
\Longrightarrow D(\overline{\lambda})=\min\left\{  D(\lambda)\mid\lambda
\in(-\Gamma_{J})\cap Y_{\operatorname{col}}^{-}\right\}  . \label{r-2b}%
\end{gather*}

\emph{(iv)} Assume that $(\overline{x},\overline{\lambda})\in\mathbb{R}%
^{n}\times\mathbb{R}^{m}$ is such that $\nabla L(\overline{x},\overline
{\lambda})=0$. Then
\begin{equation}
D(\overline{\lambda})=\left\{
\begin{array}
[c]{ll}%
\max_{\lambda\in Y_{\operatorname{col}}^{+}}D(\lambda) & \text{if }%
\overline{\lambda}\in Y_{\operatorname{col}}^{+},\\
\min_{\lambda\in Y_{\operatorname{col}}^{-}}D(\lambda) & \text{if }%
\overline{\lambda}\in Y_{\operatorname{col}}^{-}.
\end{array}
\right.  \label{r-2c}%
\end{equation}
In particular, (\ref{r-2c}) holds if $\overline{\lambda}\in Y^{+}\cup Y^{-}$
is a critical point of $D$ and $\overline{x}:=x(\overline{\lambda})$.
\end{proposition}

Proof. (i) Consider $\lambda\in Y_{\operatorname{col}}^{+};$ then there exists
$u\in\mathbb{R}^{n}$ such that $A(\lambda)u=b(\lambda)$, and so $\nabla
_{x}L(u,\lambda)=A(\lambda)u-b(\lambda)=0$. Because $L(\cdot,\lambda)$ is
convex we obtain that $L(u,\lambda)\leq L(u^{\prime},\lambda)$ for every
$u^{\prime}\in\mathbb{R}^{n}$, whence $D(\lambda)=L(u,\lambda)=\min
_{u^{\prime}\in\mathbb{R}^{n}}L(u^{\prime},\lambda)$. Of course, if
$\lambda\in Y^{+}$ then $L(\cdot,\lambda)$ is strictly convex and
$u=A(\lambda)^{-1}b(\lambda)$, and so $A(\lambda)^{-1}b(\lambda)$ is the
unique minimizer of $L(\cdot,\lambda)$ on $\mathbb{R}^{n}$. The case
$\overline{\lambda}\in Y^{-}$ is solved similarly.

(ii) Because $L(x,\cdot)$ is linear (hence concave and convex) for every
$x\in\mathbb{R}^{n}$, from (\ref{r-1}) we obtain that $D$ is concave and
u.s.c.~on $Y_{\operatorname{col}}^{+}$ as an infimum of concave continuous
functions. The argument is similar for the other situation.

(iii) Assume that $\overline{\lambda}\in Y_{\operatorname{col}}^{+}$ (hence
$\overline{\lambda}\in\Gamma_{J}\cap Y_{\operatorname{col}}^{+}$), and take
$\lambda\in\Gamma_{J}\cap Y_{\operatorname{col}}^{+}$. Using (\ref{r-1}) and
the fact that $\overline{x}\in X_{J}$, we have that
\[
D(\lambda)\leq L(\overline{x},\lambda)=q_{0}(\overline{x})+\sum_{j\in J^{c}%
}\lambda_{j}q_{j}(\overline{x})\leq q_{0}(\overline{x})=q_{0}(\overline
{x})+\left\langle \overline{\lambda},\nabla_{\lambda}L(\overline{x}%
,\overline{\lambda})\right\rangle =L(\overline{x},\overline{\lambda
})=D(\overline{\lambda}),
\]
and so $D(\overline{\lambda})=\sup_{\lambda\in\Gamma_{J}\cap
Y_{\operatorname{col}}^{+}}D(\lambda)$. The proof for $\overline{\lambda}%
\in(-\Gamma_{J})\cap Y_{\operatorname{col}}^{-}$ is similar.

(iv) One applies (iii) for $J:=\overline{1,m}$. \hfill$\square$

\medskip

Observe that $D$ is a $C^{\infty}$ function on the open set $Y$ (assumed to be
nonempty). Indeed, the operator $\varphi:\{U\in\mathfrak{M}_{n}\mid U$
invertible$\}\rightarrow\mathfrak{M}_{n}$ defined by $\varphi(U)=U^{-1}$ is
Fr\'{e}chet differentiable and $d\varphi(U)(S)=-U^{-1}SU^{-1}$ for
$U,S\in\mathfrak{M}_{n}$ with $U$ invertible. It follows that
\begin{align}
\frac{\partial D(\lambda)}{\partial\lambda_{j}} &  =\tfrac{1}{2}\left\langle
b(\lambda),A(\lambda)^{-1}A_{j}A(\lambda)^{-1}b(\lambda)\right\rangle
-\left\langle b_{j},A(\lambda)^{-1}b(\lambda)\right\rangle +c_{j}\nonumber\\
&  =\tfrac{1}{2}\left\langle x(\lambda),A_{j}x(\lambda)\right\rangle
-\left\langle b_{j},x(\lambda)\right\rangle +c_{j}=q_{j}\left(  x(\lambda
)\right)  \quad\forall j\in\overline{1,m}\label{r-d1pdq}%
\end{align}
for $\lambda\in Y_{0}$, where%
\[
x(\lambda):=A(\lambda)^{-1}b(\lambda)\quad\left(  \lambda\in Y_{0}\right)
;\label{r-xl}%
\]
hence
\begin{equation}
\nabla D(\lambda^{\prime})=\nabla_{\lambda}L(x(\lambda^{\prime}),\lambda
^{\prime})\quad\forall\lambda^{\prime}\in Y_{0}.\label{r-grpdl}%
\end{equation}
Consequently,%
\begin{equation}
\forall\lambda^{\prime}\in Y_{0}:\left[  \nabla D(\lambda^{\prime}%
)=0\iff\nabla_{\lambda}L\left(  x(\lambda^{\prime}),\lambda^{\prime}\right)
=0\iff\nabla L\left(  x(\lambda^{\prime}),\lambda^{\prime}\right)  =0\right]
.\label{r-cppdL}%
\end{equation}

A similar computation gives%
\begin{align*}
\frac{\partial^{2}D(\lambda)}{\partial\lambda_{j}\partial\lambda_{k}}= &
-\left\langle A_{j}A(\lambda)^{-1}b(\lambda),A(\lambda)^{-1}A_{k}%
A(\lambda)^{-1}b(\lambda)\right\rangle \\
&  +\left\langle A_{j}A(\lambda)^{-1}b_{k}+A_{k}A(\lambda)^{-1}b_{j}%
,A(\lambda)^{-1}b(\lambda)\right\rangle -\left\langle b_{j},A(\lambda
)^{-1}b_{k}\right\rangle \\
&  =-\left\langle A_{j}x(\lambda)-b_{j},A(\lambda)^{-1}\left(  A_{k}%
x(\lambda)-b_{k}\right)  \right\rangle \quad\forall j,k\in\overline{1,m}%
\end{align*}
for $\lambda\in Y_{0}$. Omitting $\lambda$ $(\in Y_{0})$, for $v\in
\mathbb{R}^{m}$ and $A_{v}:=\sum_{j=1}^{m}v_{j}A_{j}$, $b_{v}:=\sum_{j=1}%
^{m}v_{j}b_{j}$, we get
\[
\big\langle v,\nabla^{2}Dv\big\rangle=\sum\nolimits_{j,k=1}^{m}\frac
{\partial^{2}D}{\partial\lambda_{j}\partial\lambda_{k}}v_{j}v_{k}%
=-\left\langle A_{v}A^{-1}b-b_{v},A^{-1}\left(  A_{v}A^{-1}b-b_{v}\right)
\right\rangle .
\]
Therefore, $\nabla^{2}D(\lambda)\preceq0$ if $\lambda\in Y^{+}$ and
$\nabla^{2}D(\lambda)\succeq0$ if $\lambda\in Y^{-}$, confirming that $D$ is
concave on $Y^{+}$ and convex on $Y^{-}$.

\section{Quadratic minimization problems with equality constraints}

As mentioned above, for $J:=\overline{1,m}$, $(P_{J})$ becomes the quadratic
minimization problem

\medskip

$(P_{e})$ $~~\min$ $q_{0}(x)$ ~s.t. $x\in X_{e}:=X_{\overline{1,m}}$. \medskip

Using the previous facts we are in a position to state and prove the following result.

\begin{proposition}
\label{p1}Let $(\overline{x},\overline{\lambda})\in\mathbb{R}^{n}%
\times\mathbb{R}^{m}$.

\emph{(i)} Assume that $(\overline{x},\overline{\lambda})$ is a critical point
of $L$. Then $\overline{x}\in X_{e}$, $\overline{\lambda}\in
Y_{\operatorname{col}}$, and (\ref{r-qlpd}) holds; moreover, for
$\overline{\lambda}\in Y_{\operatorname{col}}^{+}$ we have that
\begin{equation}
q_{0}(\overline{x})=\inf_{x\in X_{e}}q_{0}(x)=L(\overline{x},\overline
{\lambda})=\sup_{\lambda\in Y_{\operatorname{col}}^{+}}D(\lambda
)=D(\overline{\lambda}), \label{r-minmaxqe}%
\end{equation}
while for $\overline{\lambda}\in Y_{\operatorname{col}}^{-}$ we have that
\begin{equation}
q_{0}(\overline{x})=\sup_{x\in X_{e}}q_{0}(x)=L(\overline{x},\overline
{\lambda})=\inf_{\lambda\in Y_{\operatorname{col}}^{-}}D(\lambda
)=D(\overline{\lambda}). \label{r-maxminqe}%
\end{equation}

\emph{(ii)} Assume that $(\overline{x},\overline{\lambda})$ is a critical
point of $L$ with $\overline{\lambda}\in Y_{0}$. Then $\nabla D(\overline
{\lambda})=0$ and $\overline{x}=A(\overline{\lambda})^{-1}b(\overline{\lambda
})$; moreover, $\overline{x}$ is the unique global minimizer of $q_{0}$ on
$X_{e}$ when $\overline{\lambda}\in Y^{+}$, and $\overline{x}$ is the unique
global maximizer of $q_{0}$ on $X_{e}$ when $\overline{\lambda}\in Y^{-}$.

Conversely, assume that $\overline{\lambda}\in Y_{0}$ is a critical point of
$D$. Then $(\overline{x},\overline{\lambda})$ is a critical point of $L$,
where $\overline{x}=A(\overline{\lambda})^{-1}b(\overline{\lambda})$;
consequently \emph{(i)} and \emph{(ii)} apply.
\end{proposition}

Proof. (i) Assume that $(\overline{x},\overline{\lambda})$ is a critical point
of $L$; hence $\nabla_{x}L(\overline{x},\overline{\lambda})=0$ and
$\nabla_{\lambda}L(\overline{x},\overline{\lambda})=0$. Using Lemma
\ref{lem-qperfdual} we obtain that $\overline{\lambda}\in
Y_{\operatorname{col}}$, $\overline{x}\in X_{e}$, and (\ref{r-qlpd}) holds.

Assume moreover that $\overline{\lambda}\in Y_{\operatorname{col}}^{+}$.
Because $L(\cdot,\overline{\lambda})$ is convex, its infimum is attained at
$\overline{x}$. Therefore, for $x\in X_{e}$ we have that $q_{0}(\overline
{x})=L(\overline{x},\overline{\lambda})\leq L(x,\overline{\lambda})=q_{0}(x)$,
and so $q_{0}(\overline{x})=\inf_{x\in X_{e}}q_{0}(x)$. Using Proposition
\ref{lem-pd}~(iii) for $J:=\overline{1,m}$ (hence $\Gamma_{J}=\mathbb{R}^{m}%
$), we get the last equality in (\ref{r-minmaxqe}). Hence (\ref{r-minmaxqe}) holds.

The proof of (\ref{r-maxminqe}) in the case $\overline{\lambda}\in
Y_{\operatorname{col}}^{-}$ is similar; an alternative proof is to apply the
previous case for $q_{j}$ replaced by $-q_{j}$ and $\overline{\lambda}_{j}$ by
$-\overline{\lambda}_{j}$ for $j\in\overline{1,m}.$

(ii) Assume that $(\overline{x},\overline{\lambda})$ is a critical point of
$L$ with $\overline{\lambda}\in Y_{0}$. Since $A(\overline{\lambda}%
)\overline{x}-b(\overline{\lambda})=\nabla_{x}L(\overline{x},\overline
{\lambda})=0$, clearly $\overline{x}=x(\overline{\lambda})=A(\overline
{\lambda})^{-1}b(\overline{\lambda})$. Using (\ref{r-grpdl}) we obtain that
$\nabla D(\overline{\lambda})=\nabla_{\lambda}L(\overline{x},\overline
{\lambda})=0$.

Moreover, suppose that $\overline{\lambda}\in Y^{+}$. Then $L(\cdot
,\overline{\lambda})$ is strictly convex, and so $q_{0}(\overline
{x})=L(\overline{x},\overline{\lambda})<L(x,\overline{\lambda})=q_{0}(x)$ for
$x\in X_{e}\setminus\{\overline{x}\}$. Hence $\overline{x}$ is the unique
global minimizer of $q_{0}$ on $X_{e}$. The proof in the case $\overline
{\lambda}\in Y^{-}$ is similar.

Conversely, let $\overline{\lambda}\in Y_{0}$ be a critical point of $D$ and
take $\overline{x}:=A(\overline{\lambda})^{-1}b(\overline{\lambda});$ then
$\nabla_{x}L(\overline{x},\overline{\lambda})=0$ by (\ref{r-nxL}). Using
(\ref{r-d1pdq}) we obtain that $\overline{x}\in X_{e}$, and so $\nabla
_{\lambda}L(\overline{x},\overline{\lambda})=0$. Therefore, $(\overline
{x},\overline{\lambda})$ is a critical point of $L$. \hfill$\square$

\medskip

The next example shows that $(P_{e})$ might have several solutions when
$\overline{\lambda}\in Y_{\operatorname{col}}^{+}$.

\begin{example}
\label{ex-qe}Take $q_{0}(x,y):=xy$, $q_{1}(x,y):=\tfrac{1}{2}(x^{2}+y^{2}-1)$
for $x,y\in\mathbb{R}$. Then $L(x,y,\lambda)=xy+\tfrac{\lambda}{2}\left(
x^{2}+y^{2}-1\right)  $. It follows that $A(\lambda)=\left(
\begin{array}
[c]{ll}%
\lambda & 1\\
1 & \lambda
\end{array}
\right)  $, $b(\lambda)=\left(
\begin{array}
[c]{l}%
0\\
0
\end{array}
\right)  $, $c(\lambda)=-\tfrac{1}{2}\lambda$, $Y_{0}=\mathbb{R}%
\setminus\{-1,1\}$, $Y^{+}=-Y^{-}=(1,\infty)$, $Y_{\operatorname{col}%
}=\mathbb{R}$, $Y_{\operatorname{col}}^{+}=-Y_{\operatorname{col}}%
^{-}=[1,\infty)$, $D(\lambda)=-\tfrac{1}{2}\lambda$. Clearly, $D$ has not
critical points, and the only critical points of $L$ are $(\pm2^{-1/2}%
,\mp2^{-1/2},1)$ and $(\pm2^{-1/2},\pm2^{-1/2},-1)$. For $(\pm2^{-1/2}%
,\mp2^{-1/2},1)$ we can apply Proposition \ref{p1}~(i) with $\overline
{\lambda}:=1\in Y_{\operatorname{col}}^{+}$, and so both $\pm2^{-1/2}(1,-1)$
are solutions for problem $(P_{e})$, while for $(\pm2^{-1/2},\pm2^{-1/2},-1)$
we can apply Proposition \ref{p1}~(i) with $\overline{\lambda}:=-1\in
Y_{\operatorname{col}}^{-}$, and so $\pm2^{-1/2}(1,1)$ are global maximizers
of $q_{0}$ on $X_{e}$.
\end{example}

\section{Quadratic minimization problems with equality and inequality
constraints}

Let us consider now the general quadratic minimization problem $(P_{J})$
considered at the beginning of Section \ref{sec1}. To $(P_{J})$ we associate
the sets
\begin{gather*}
Y^{J}:=\Gamma_{J}\cap Y_{0},\quad Y^{J+}:=\Gamma_{J}\cap Y^{+},\quad
Y^{J-}:=(-\Gamma_{J})\cap Y^{-},\\
Y_{\operatorname{col}}^{J}:=\Gamma_{J}\cap Y_{\operatorname{col}},\quad
Y_{\operatorname{col}}^{J+}:=\Gamma_{J}\cap Y_{\operatorname{col}}^{+},\quad
Y_{\operatorname{col}}^{J-}:=(-\Gamma_{J})\cap Y_{\operatorname{col}}^{-},
\end{gather*}
where $Y_{0}$, $Y^{+}$ and $Y^{-}$, $Y_{\operatorname{col}}$,
$Y_{\operatorname{col}}^{+}$ and $Y_{\operatorname{col}}^{-}$, are defined in
(\ref{r-s0}), (\ref{r-s0pm}), (\ref{r-yext}) and (\ref{r-yextp}),
respectively. Unlike $Y_{0}$, $Y^{+}$, $Y^{-}$, the sets $Y^{J}$, $Y^{J+}$ and
$Y^{J-}$ are (generally) not open. Because $Y^{+}$, $Y_{\operatorname{col}%
}^{+}$ and $Y_{\operatorname{col}}^{-}$ are convex, so are $Y^{J+}$,
$Y_{\operatorname{col}}^{J+}$ and $Y_{\operatorname{col}}^{J-}$, and so
$L(\cdot,\lambda)$ is (strictly) convex on $Y_{\operatorname{col}}^{J+}$
$(Y^{J+})$ and (strictly) concave on $Y_{\operatorname{col}}^{J-}$ $(Y^{J-});$
moreover, $\operatorname*{int}Y^{J+}=\operatorname*{int}Y_{\operatorname{col}%
}^{J+}$ ($\operatorname*{int}Y^{J-}=\operatorname*{int}Y_{\operatorname{col}%
}^{J-}$) provided $Y^{J+}\neq\emptyset$ ($\operatorname*{int}Y^{J-}%
\neq\emptyset$).

As observed already, for $J=\overline{1,m}$ we have that $\Gamma
_{J}=\mathbb{R}^{m}$, and so $Y^{J}$, $Y^{J+}$, $Y^{J-}$,
$Y_{\operatorname{col}}^{J}$, $Y_{\operatorname{col}}^{J+}$ and
$Y_{\operatorname{col}}^{J-}$ reduce to $Y_{0}$, $Y^{+}$, $Y^{-}$,
$Y_{\operatorname{col}}$, $Y_{\operatorname{col}}^{+}$ and
$Y_{\operatorname{col}}^{-}$, respectively.

Suggested by the well known necessary optimality conditions for minimization
problems with equality and inequality constraints, we say that $(\overline
{x},\overline{\lambda})\in\mathbb{R}^{n}\times\mathbb{R}^{m}$ is a $J$-LKKT
point of $L$ (that is a Lagrange--Karush--Kuhn--Tucker\footnote{It seems that
the term Lagrange--Karush--Kuhn--Tucker multiplier was introduced by J.-P.
Penot in \cite{Pen:81}.} point of $L$ with respect to $J$) if $\nabla
_{x}L(\overline{x},\overline{\lambda})=0$ and
\[
\textstyle\left[  \forall j\in J^{c}:\overline{\lambda}_{j}\geq0~~\wedge
~~\frac{\partial L}{\partial\lambda_{j}}(\overline{x},\overline{\lambda}%
)\leq0~~\wedge~~\overline{\lambda}_{j}\cdot\frac{\partial L}{\partial
\lambda_{j}}(\overline{x},\overline{\lambda})=0\right]  ~\wedge~\left[
\forall j\in J:\frac{\partial L}{\partial\lambda_{j}}(\overline{x}%
,\overline{\lambda})=0\right]  ,\label{r-kkt-lqm}%
\]
or, equivalently,
\begin{equation}
\overline{x}\in X_{J}~~\wedge~~\overline{\lambda}\in\Gamma_{J}~~\wedge
~~\left[  \forall j\in J^{c}:\overline{\lambda}_{j}q_{j}(\overline
{x})=0\right]  ;\label{r-kkt-pqei}%
\end{equation}
we say that $\overline{x}\in\mathbb{R}^{n}$ is a $J$-LKKT point for $(P_{J})$
if there exists $\overline{\lambda}\in\mathbb{R}^{m}$ such that $(\overline
{x},\overline{\lambda})$ verifies (\ref{r-kkt-pqei}); moreover, for $D$
defined in (\ref{r-pd0}), we say that $\overline{\lambda}\in Y_{0}$ is a
$J$-LKKT point for $D$ if%
\begin{equation}
\textstyle\left[  \forall j\in J^{c}:\overline{\lambda}_{j}\geq0~~\wedge
~~\frac{\partial D}{\partial\lambda_{j}}(\overline{\lambda})\leq
0~~\wedge~~\overline{\lambda}_{j}\cdot\frac{\partial D}{\partial\lambda_{j}%
}(\overline{\lambda})=0\right]  ~\wedge~\left[  \forall j\in J:\frac{\partial
D}{\partial\lambda_{j}}(\overline{\lambda})=0\right]  .\label{r-kkt-dqm}%
\end{equation}

Of course, when $J=\overline{1,m}$, $(\overline{x},\overline{\lambda}%
)\in\mathbb{R}^{n}\times\mathbb{R}^{m}$ is a $J$-LKKT point of $L$ iff $\nabla
L(\overline{x},\overline{\lambda})=0$, while $\overline{\lambda}\in Y_{0}$ is
a $J$-LKKT point for $D$ iff $\nabla D(\overline{\lambda})=0$.

\begin{remark}
\label{rem-kktLkktpd}Notice that $\overline{\lambda}\in Y_{0}$ is a $J$-LKKT
point of $D$ if and only if $(x(\overline{\lambda}),\overline{\lambda})$ is a
$J$-LKKT point of $L$; for this just take into account (\ref{r-d1pdq}).
Moreover, taking into account (\ref{r-cppdL}), if $\overline{\lambda}\in
Y_{0}$ is a critical point of $D$ then $\overline{\lambda}$ is a $J$-LKKT
point of $D$ and $(x(\overline{\lambda}),\overline{\lambda})$ is a $J$-LKKT
point of $L$ (being a critical point of $L$).
\end{remark}

In general, for distinct $J$ and $J^{\prime}$, the sets of $J$-LKKT and
$J^{\prime}$-LKKT points of $L$ (resp. $D$) are not comparable. For comparable
$J$ and $J^{\prime}$ we have the following result whose simple proof is
omitted; its second part follows from the first one and the previous remark.

\begin{lemma}
\label{fact1}Let $J\subset J^{\prime}\subset\overline{1,m}$ and $(\overline
{x},\overline{\lambda})\in\mathbb{R}^{n}\times\mathbb{R}^{m}$

\emph{(i)} If $(\overline{x},\overline{\lambda})$ is a $J^{\prime}$-LKKT point
of $L$ and $\overline{\lambda}_{j}\geq0$ for all $j\in J^{\prime}\setminus J$,
then $(\overline{x},\overline{\lambda})$ is a $J$-LKKT point of $L$.
Conversely, if $(\overline{x},\overline{\lambda})$ is a $J$-LKKT point of $L$
and $\overline{\lambda}_{j}>0$ for all $j\in J^{\prime}\setminus J$, then
$(\overline{x},\overline{\lambda})$ is a $J^{\prime}$-LKKT point of $L$.

\emph{(ii)} If $\overline{\lambda}\in Y_{0}$ is a $J^{\prime}$-LKKT point of
$D$ and $\overline{\lambda}_{j}\geq0$ for all $j\in J^{\prime}\setminus J$,
then $\overline{\lambda}$ is a $J$-LKKT point of $D$. Conversely, if
$\overline{\lambda}$ is a $J$-LKKT point of $D$ and and $\overline{\lambda
}_{j}>0$ for all $j\in J^{\prime}\setminus J$, then $\overline{\lambda}$ is a
$J^{\prime}$-LKKT point of $D$.
\end{lemma}

The result below corresponds to Proposition \ref{p1}.

\begin{proposition}
\label{p1ei}Let $(\overline{x},\overline{\lambda})\in\mathbb{R}^{n}%
\times\mathbb{R}^{m}$.

\emph{(i)} Assume that $(\overline{x},\overline{\lambda})$ is a $J$-LKKT point
of $L$. Then $\overline{x}$ is a $J$-LKKT point of $(P_{J})$, $\overline{x}\in
X_{J}$, $\overline{\lambda}\in Y_{\operatorname{col}}^{J}$, and (\ref{r-qlpd})
holds; moreover, if $\overline{\lambda}\in Y_{\operatorname{col}}^{J+}$ then
\begin{equation}
q_{0}(\overline{x})=\inf_{x\in X_{J}}q_{0}(x)=L(\overline{x},\overline
{\lambda})=\sup_{\lambda\in Y_{\operatorname{col}}^{J+}}D(\lambda
)=D(\overline{\lambda}). \label{r-minmaxqei}%
\end{equation}

\emph{(ii) }Assume that $(\overline{x},\overline{\lambda})$ is a $J$-LKKT
point of $L$ with $\overline{\lambda}\in Y_{0}$ (or, equivalently,
$\overline{\lambda}\in Y^{J}$). Then $\overline{x}=A(\overline{\lambda}%
)^{-1}b(\overline{\lambda})$, and $\overline{\lambda}$ is a $J$-LKKT point of
$D$; moreover, $\overline{x}$ is the unique global minimizer of $q_{0}$ on
$X_{J}$ if $\overline{\lambda}\in Y^{J+}$.

Conversely, assume that $\overline{\lambda}\in Y_{0}$ is a $J$-LKKT point of
$D$. Then $(\overline{x},\overline{\lambda})$ is a $J$-LKKT point of $L$,
where $\overline{x}:=A(\overline{\lambda})^{-1}b(\overline{\lambda})$.
Consequently, \emph{(i)} and \emph{(ii)} apply.

\emph{(iii)} Assume that $\overline{\lambda}\in Y^{J+}$. Then%
\[
D(\overline{\lambda})=\sup_{\lambda\in Y_{\operatorname{col}}^{J+}}%
D(\lambda)\Longleftrightarrow D(\overline{\lambda})=\sup_{\lambda\in Y^{J+}%
}D(\lambda)\Longleftrightarrow\overline{\lambda}\text{ is a $J$-LKKT point of
$D$.}%
\]

\end{proposition}

Proof. (i) By hypothesis, (\ref{r-kkt-pqei}) holds. The fact that
$\overline{x}$ is a $J$-LKKT point of $(P_{J})$ is obvious from its very
definition; hence $\overline{x}\in X_{J}$. On the other hand, because
$(\overline{x},\overline{\lambda})$ is a $J$-LKKT point of $L$ we have that
$\overline{\lambda}\in Y_{\operatorname{col}}^{J}$ and (\ref{r-qlpd}) holds by
Lemma \ref{lem-qperfdual}.

Assume that $\overline{\lambda}\in Y_{\operatorname{col}}^{J+}$ $(=\Gamma
_{J}\cap Y_{\operatorname{col}}^{+})$. The last equality in (\ref{r-minmaxqei}%
) follows from Proposition \ref{lem-pd}~(iii). Because $L(\cdot,\overline
{\lambda})$ is convex, its infimum is attained at $\overline{x}$. Therefore,
for $x\in X_{J}$ we have that
\[
q_{0}(\overline{x})=L(\overline{x},\overline{\lambda})\leq L(x,\overline
{\lambda})=q_{0}(x)+\sum\nolimits_{j=1}^{m}\overline{\lambda}_{j}q_{j}(x)\leq
q_{0}(x), \label{r-siqei}%
\]
whence $q_{0}(\overline{x})=\inf_{x\in X_{i}}q_{0}(x)$. Hence
(\ref{r-minmaxqei}) holds.

(ii) Because $(\overline{x},\overline{\lambda})$ is a $J$-LKKT point of $L$
with $\overline{\lambda}\in Y_{0}$, we have that $A(\overline{\lambda
})\overline{x}-b(\overline{\lambda})=\nabla_{x}L(\overline{x},\overline
{\lambda})=0$, and so $\overline{x}=x(\overline{\lambda})$. As observed in
Remark \ref{rem-kktLkktpd}, (\ref{r-kkt-dqm}) is verified.

Suppose now that moreover that $\overline{\lambda}\in Y^{+}$ (and so Then
$L(\cdot,\overline{\lambda})$ is strictly convex, and so $q_{0}(\overline
{x})=L(\overline{x},\overline{\lambda})<L(x,\overline{\lambda})\leq q_{0}(x)$
for $x\in X_{J}\setminus\{\overline{x}\}$. Hence $\overline{x}$ is the unique
global minimizer of $q_{0}$ on $X_{J}$.

Conversely, let $\overline{\lambda}\in Y_{0}$ be a $J$-LKKT point of $D$, and
take $\overline{x}:=x(\overline{\lambda});$ then $(\overline{x},\overline
{\lambda})$ is a $J$-LKKT point of $L$ by Remark \ref{rem-kktLkktpd}.

(iii) If $\overline{\lambda}$ is a $J$-LKKT point of $D$, we have that
$D(\overline{\lambda})=\sup_{\lambda\in Y_{\operatorname{col}}^{J+}}%
D(\lambda)$ by Remark \ref{rem-kktLkktpd} and (i), while $D(\overline{\lambda
})=\sup_{\lambda\in Y_{\operatorname{col}}^{J+}}D(\lambda)$ implies
$D(\overline{\lambda})=\sup_{\lambda\in Y^{J+}}D(\lambda)$ because
$Y^{J+}\subset Y_{\operatorname{col}}^{J+}$. Assume that $D(\overline{\lambda
})=\sup_{\lambda\in Y^{J+}}D(\lambda)$. Setting $Q:=-D$, we have that $Q$ is
convex and $\overline{\lambda}$ is a global minimizer of $Q$ on (the convex
set) $Y^{J+}$. Using \cite[Prop.\ 4]{Zal:89} we have that
\[
0\leq Q^{\prime}(\overline{\lambda},\lambda-\overline{\lambda}):=\lim
_{t\rightarrow0+}\frac{Q(\overline{\lambda}+t(\lambda-\overline{\lambda
}))-Q(\overline{\lambda})}{t}=\left\langle \lambda-\overline{\lambda},\nabla
Q(\overline{\lambda})\right\rangle \quad\forall\lambda\in Y^{J+}.
\]
It follows that $\left\langle y,v\right\rangle \leq0$ for all $y\in
\mathbb{R}_{+}(Y^{J+}-\overline{\lambda})$, where $v:=\nabla D(\overline
{\lambda})$. Because $\Gamma_{J}$ and $Y^{+}$ are convex sets, $Y^{J+}%
=\Gamma_{J}\cap Y^{+}$, and $\overline{\lambda}\in\operatorname*{int}%
Y^{+}=Y^{+}$, we have that
\begin{align*}
\mathbb{R}_{+}(Y^{J+}-\overline{\lambda})  &  =\mathbb{R}_{+}\left[
(\Gamma_{J}-\overline{\lambda})\cap(Y^{+}-\overline{\lambda})\right]
=\mathbb{R}_{+}(\Gamma_{J}-\overline{\lambda})\\
&  =\left\{  \mu\in\mathbb{R}^{m}\mid\forall j\in J^{c}:\overline{\lambda}%
_{j}=0\Rightarrow\mu_{j}\geq0\right\}  .
\end{align*}
Therefore, $\frac{\partial D}{\partial\lambda_{j}}(\overline{\lambda}%
)=v_{j}=0$ for $j\in J\cup\{j\in J^{c}\mid\overline{\lambda}_{j}>0\}$ and
$\frac{\partial D}{\partial\lambda_{j}}(\overline{\lambda})=v_{j}\leq0$ for
$j\in\{j^{\prime}\in J\mid\overline{\lambda}_{j^{\prime}}=0\}$. This shows
that condition (\ref{r-kkt-dqm}) is verified. \hfill$\square$

\begin{corollary}
\label{fact2}Let $\emptyset\neq J\subset\overline{1,m}$ and let $(\overline
{x},\overline{\lambda})\in\mathbb{R}^{n}\times\mathbb{R}^{m}$ be a $J$-LKKT
point of $L$ such that $A(\overline{\lambda})\succeq0;$ hence $\overline{x}$
$\in X_{J}$, $\overline{\lambda}\in Y_{\operatorname{col}}^{J+}$ and
(\ref{r-minmaxqei}) holds. If $J_{\geq}:=\{j\in J\mid\overline{\lambda}%
_{j}\geq0\}$ is nonempty, then $(\overline{x},\overline{\lambda})$ is a
$(J\setminus J_{\geq})$-LKKT point of $L$, and so $\overline{x}$ is a global
minimizer of $q_{0}$ on $X_{J\setminus J_{\geq}}\supset X_{J}$.
\end{corollary}

Proof. The first assertion holds by Proposition \ref{p1ei}~(i) because
$\overline{\lambda}\in Y_{\operatorname{col}}^{J+}$. In what concerns the
second assertion, it is sufficient to observe that for $j\in J^{c}\cup
J_{\geq}=(J\setminus J_{\geq})^{c}$ we have that $\overline{\lambda}_{j}\geq
0$, and $\overline{\lambda}_{j}\cdot\frac{\partial L}{\partial\lambda_{j}%
}(\overline{x},\overline{\lambda})=0$ by the definition of a $J$-LKKT point of
$L$, then to apply Proposition \ref{p1ei}~(i) for $J$ replaced by $J\setminus
J_{\geq}$. \hfill$\square$

\begin{corollary}
\label{fact3}If $(\overline{x},\overline{\lambda})\in\mathbb{R}^{n}%
\times\mathbb{R}^{m}$ is a critical point of $L$ (in particular if
$\overline{\lambda}\in Y_{0}$ is a critical point of $D$ and $\overline
{x}:=x(\overline{\lambda})$), then $(\overline{x},\overline{\lambda})$ is a
$J^{c}$-LKKT point of $L$, where $J:=\{j\in\overline{1,m}\mid\overline
{\lambda}_{j}\geq0\}$. Consequently, if moreover $A(\overline{\lambda})\geq0$,
then $\overline{x}$ $(\in X_{e})$ is a global minimizer of $q_{0}$ on
$X_{J^{c}}\supset X_{e}$.
\end{corollary}

Proof. Apply Corollary \ref{fact2} for $J:=\overline{1,m}$. \hfill$\square$

\medskip

The next result is the variant of Proposition \ref{p1ei} for maximizing
$q_{0}$ on $X_{J}$.

\begin{proposition}
\label{p1eimax}Let $(\overline{x},\overline{\lambda})\in\mathbb{R}^{n}%
\times\mathbb{R}^{m}$.

\emph{(i)} Assume that $\nabla_{x}L(\overline{x},\overline{\lambda})=0$ and
the condition%
\begin{equation}
\textstyle \left[  \forall j\in J^{c}:\overline{\lambda}_{j}\leq
0~~\wedge~~\frac{\partial L}{\partial\lambda_{j}}(\overline{x},\overline
{\lambda})\leq0~~\wedge~~\overline{\lambda}_{j}\cdot\frac{\partial L}%
{\partial\lambda_{j}}(\overline{x},\overline{\lambda})=0\right]
~\wedge~\left[  \forall j\in J:\frac{\partial L}{\partial\lambda_{j}%
}(\overline{x},\overline{\lambda})=0\right]  \label{r-lkktmax}%
\end{equation}
is verified. Then $\overline{x}$ $\in X_{J}$, $\overline{\lambda}\in
Y_{\operatorname{col}}$, and
\begin{equation}
\textstyle \left[  \forall j\in J^{c}:\overline{\lambda}_{j}\leq
0~~\wedge~~\frac{\partial D}{\partial\lambda_{j}}(\overline{\lambda}%
)\leq0~~\wedge~~\overline{\lambda}_{j}\cdot\frac{\partial D}{\partial
\lambda_{j}}(\overline{\lambda})=0\right]  ~\wedge~\left[  \forall j\in
J:\frac{\partial D}{\partial\lambda_{j}}(\overline{\lambda})=0\right]  ;
\label{r-qlpdmax}%
\end{equation}
moreover, if $\overline{\lambda}\in Y_{\operatorname{col}}^{-}$ (or
equivalently $\overline{\lambda}\in Y_{\operatorname{col}}^{J-}$), then
\[
q_{0}(\overline{x})=\sup_{x\in X_{J}}q_{0}(x)=L(\overline{x},\overline
{\lambda})=\inf_{\lambda\in Y_{\operatorname{col}}^{J-}}D (\lambda
)=D(\overline{\lambda}).
\]

\emph{(ii) }Assume that $\overline{\lambda}\in Y_{0}$, $\nabla_{x}%
L(\overline{x},\overline{\lambda})=0$ and $(\overline{x},\overline{\lambda})$
verifies (\ref{r-lkktmax}). Then $\overline{x}=x(\overline{\lambda})$, and
$\overline{\lambda}$ verifies condition (\ref{r-qlpdmax}); moreover,
$\overline{x}$ is the unique global maximizer of $q_{0}$ on $X_{J}$ if
$\overline{\lambda}\in Y^{J-}$.

\emph{(iii)} Assume that $\overline{\lambda}\in Y^{J-}$. Then%
\[
D(\overline{\lambda})=\inf_{\lambda\in Y^{J-}}D(\lambda)\Longleftrightarrow
D(\overline{\lambda})=\inf_{\lambda\in Y_{\operatorname{col}}^{J-}}%
D(\lambda)\Longleftrightarrow\overline{\lambda}\text{ verifies condition
(\ref{r-qlpdmax}).}%
\]

\end{proposition}

The proof of the above result is an easy adaptation of the proof of
Proposition \ref{p1ei}, so we omit it.

\section{Quadratic minimization problems with inequality constraints}

We consider now the particular case of $(P_{J})$ in which $J=\emptyset;$ the
problem is denoted by $(P_{i})$ and the set of its feasible solutions by
$X_{i}$. In this case $\Gamma_{J}=\mathbb{R}_{+}^{m}$, and the sets $Y^{J}$,
$Y^{J+}$, $Y^{J-}$, $Y_{\operatorname{col}}^{J}$, $Y_{\operatorname{col}}%
^{J+}$, $Y_{\operatorname{col}}^{J+}$ and $Y_{\operatorname{col}}^{J-}$ are
denoted by $Y^{i}$, $Y^{i+}$, $Y^{i-}$, $Y_{\operatorname{col}}^{i}$,
$Y_{\operatorname{col}}^{i+}$ and $Y_{\operatorname{col}}^{i-}$, respectively.
Moreover, in this situation we shall use KKT instead of $J$-LKKT. So, we say
that $(\overline{x},\overline{\lambda})\in\mathbb{R}^{n}\times\mathbb{R}^{m}$
is a Karush--Kuhn--Tucker point of $L$ if $\nabla_{x}L(\overline{x}%
,\overline{\lambda})=0$ and
\[
\overline{\lambda}\in\mathbb{R}_{+}^{m}~~\wedge~~\nabla_{\lambda}%
L(\overline{x},\overline{\lambda})\in\mathbb{R}_{-}^{m}~~\wedge
~~\big \langle\overline{\lambda},\nabla_{\lambda}L(\overline{x},\overline
{\lambda})\big\rangle=0,\label{r-kkt-lq}%
\]
or, equivalently,%
\begin{equation}
\overline{x}\in X_{i}~~\wedge~~\overline{\lambda}\in\mathbb{R}_{+}^{m}%
~~\wedge~~\left[  \forall j\in\overline{1,m}:\overline{\lambda}_{j}%
q_{j}(\overline{x})=0\right]  ;\label{r-kkt-pqi}%
\end{equation}
we say that $\overline{x}$ is a KKT point for $(P_{i})$ if there exists
$\overline{\lambda}\in\mathbb{R}^{m}$ such that (\ref{r-kkt-pqi}) holds; we
say that $\overline{\lambda}\in Y_{0}$ is a KKT point for $D$ if%
\[
\overline{\lambda}\in\mathbb{R}_{+}^{m}~~\wedge~~\nabla D(\overline{\lambda
})\in\mathbb{R}_{-}^{m}~~\wedge~~\big \langle\overline{\lambda},\nabla
D(\overline{\lambda})\big\rangle=0.\label{r-kkt-dqi}%
\]

Proposition \ref{p1ei} becomes the next result when $J=\emptyset$.

\begin{proposition}
\label{p1i}Let $(\overline{x},\overline{\lambda})\in\mathbb{R}^{n}%
\times\mathbb{R}^{m}$.

\emph{(i)} Assume that $(\overline{x},\overline{\lambda})$ is a KKT point of
$L$. Then $\overline{x}$ is a KKT point of $(P_{i})$, and so $\overline{x}\in
X_{i}$, $\overline{\lambda}\in Y_{\operatorname{col}}^{i}$, and (\ref{r-qlpd})
holds; moreover, for $\overline{\lambda}\in Y_{\operatorname{col}}^{i+}$ we
have that
\begin{equation}
q_{0}(\overline{x})=\inf_{x\in X_{i}}q_{0}(x)=L(\overline{x},\overline
{\lambda})=\sup_{\lambda\in Y_{\operatorname{col}}^{i+}}D(\lambda
)=D(\overline{\lambda}). \label{r-minmaxqi}%
\end{equation}

\emph{(ii)} Assume that $(\overline{x},\overline{\lambda})$ is a KKT point of
$L$ with $\overline{\lambda}\in Y_{0}$. Then $\overline{x}=x(\overline
{\lambda})$ and $\overline{\lambda}$ is a KKT point of $D$; moreover,
$\overline{x}$ is the unique global minimizer of $q_{0}$ on $X_{i}$ provided
$\overline{\lambda}\in Y^{i+}$.

Conversely, assume that $\overline{\lambda}\in Y_{0}$ is a KKT point of $D$.
Then $(\overline{x},\overline{\lambda})$ is a KKT point of $L$, where
$\overline{x}:=A(\overline{\lambda})^{-1}b(\overline{\lambda})$.

\emph{(iii)} Assume that $\overline{\lambda}\in Y^{i+}$. Then
\[
D(\overline{\lambda})=\sup_{\lambda\in Y_{\operatorname{col}}^{i+}}%
D(\lambda)\Longleftrightarrow D(\overline{\lambda})=\sup_{\lambda\in Y^{i+}%
}D(\lambda)\Longleftrightarrow\overline{\lambda}\text{ is a KKT point of $D$.}%
\]

\end{proposition}

\begin{remark}
\label{r-jrw}Jeyakumar, Rubinov and Wu (see \cite[Prop.~3.2]{JeyRubWu:07})
proved that $\overline{x}$ is a (global) solution of $(P_{i})$ when there
exists $\overline{\lambda}\in Y_{\operatorname{col}}^{i}$ is a KKT point of
$L$; this result was established previously by Hiriart-Urruty in
\cite[Th.~4.6]{Hir:98} when $m=2$.
\end{remark}

\begin{remark}
\label{rem-adv}Having in view Propositions \ref{p1}, \ref{p1ei}, \ref{p1i}, it
is more advantageous to use their versions (i) than the second part of (ii)
with $\overline{\lambda}\in Y_{0}$ because in versions (i) one must know only
the Lagrangian (hence only the data of the problems), and this provides both
$\overline{x}$ and $\overline{\lambda}$, without needing to calculate
effectively $D$, then to determine $\overline{\lambda}$ (and after that,
$\overline{x}$). Using $D$ could be useful, maybe, if the number of
constraints is much smaller than $n$. As seen in the proofs, the consideration
of the dual function is not essential in finding the optimal solutions of the
primal problem(s).
\end{remark}

\section{Comparisons with results on quadratic optimization problems obtained
by using CDT\label{sec4}}

In this section we analyze results obtained by DY Gao and his collaborators in
papers dedicated to quadratic optimization problems, or as particular cases of
more general results. The main tool to identify the papers where quadratic
problems are considered was to look in the survey papers like \cite{Gao:03a},
\cite{Gao:09} (which is almost the same as \cite{Gao:08}, both of them being
cited in Gao's papers), \cite{GaoShe:09} (which is very similar to
\cite{Gao:09b}), as well as in the recent book \cite{GaoLatRua:17}.

We present the results in chronological order using our notations (when
possible) and with equivalent formulations; however, sometimes we quote the
original formulations to feel also the flavor of those papers. When we have
not notations for some sets we introduce them, often as in the respective
papers; similarly for some notions. Because $c_{0}$ in the definition of
$q_{0}$ may be taken always to be $0$, we shall not mention it in the sequel.

Before beginning our analysis we consider it is worth having in view the
following remark from the very recent paper \cite{RuaGao:17b} and to observe
that there is not an assumption that some multiplier $\overline{\lambda}_{j}$
be non null in Propositions \ref{p1}, \ref{p1ei} and \ref{p1i}.

\begin{quote}
``\emph{Remark 1}. As we have demonstrated that by the generalized canonical
duality (32), all KKT conditions can be recovered for both equality and
inequality constraints. Generally speaking, the nonzero Lagrange multiplier
condition for the linear equality constraint is usually ignored in
optimization textbooks. But it can not be ignored for nonlinear constraints.
It is proved recently [26] that the popular augmented Lagrange multiplier
method can be used mainly for linear constrained problems. Since the
inequality constraint $\mu\not =0$ produces a nonconvex feasible set
$\mathcal{E}_{a}^{\ast}$, this constraint can be replaced by either $\mu<0$ or
$\mu>0$. But the condition $\mu<0$ is corresponding to $y\circ(y-e_{K})\geq0$,
this leads to a nonconvex open feasible set for the primal problem. By the
fact that the integer constraints $y_{i}(y_{i}-1)=0$ are actually a special
case (boundary) of the boxed constraints $0\leq y_{i}\leq1$, which is
corresponding to $y\circ(y-e_{K})\geq0$, we should have $\mu>0$ (see [8] and
[12, 16]). In this case, the KKT condition (43) should be replaced by

$\mu>0,~~y\circ(y-e_{K})\leq0,~~\mu^{T}[y\circ(y-e_{K})]=0.\quad$ (47)

\noindent Therefore, as long as $\mu\neq0$ is satisfied, the complementarity
condition in (47) leads to the integer condition $y\circ(y-e_{K})=0$.
Similarly, the inequality $\tau\neq0$ can be replaced by $\tau>0$."
\end{quote}


Notice that many papers (co-) authored by DY Gao, mostly in those made public
in the last five years, the multipliers corresponding to nonlinear constraints
(but not only) are assumed to be positive. So, in most cases Eq.~(\ref{r-2c})
is true. Moreover, it is worth observing that $\overline{x}\in X_{J}$ is a
local minimizer as well as a local maximizer of $q_{0}$ on $X_{J}$ whenever
$X_{J}$ is a finite set; this is the case in many optimization problems
mentioned in this section.

\bigskip

The quadratic problem considered by Gao in \cite[Sect.~5.1]{Gao:03a} is of
type $(P_{i})$ in which $A_{1}:=I_{n}:=\operatorname*{diag}e$ with
$e:=(1,...,1)^{T}\in\mathbb{R}^{n}$, $b_{1}=0$, $c_{1}<0$, $A_{j}=0$ for
$j\in\overline{2,m}$. Below, $X_{i1}:=\{x\in X_{i}\mid q_{1}(x)=0\}$ and
$Y_{1}^{i+}:=\{\lambda\in Y^{i+}\mid\lambda_{1}>0\}$.

\smallskip Theorem 4 in \cite{Gao:03a} (attributed to \cite{Gao:03b}) asserts:
\emph{Let $\overline{\lambda}\in Y^{i}$ be a KKT point of $D$ and
$\overline{x}:=x(\overline{\lambda})$. Then $\overline{x}$ is a KKT point of
$(P_{i})$ and $q_{0}(\overline{x})=D(\overline{\lambda})$.}

\smallskip Theorem 6 in \cite{Gao:03a} asserts: \emph{Assume that $A_{0}$ has
at least one negative eigenvalue and $(\overline{x},\overline{\lambda})$ is a
KKT point of $L$. If $\overline{\lambda}\in Y_{1}^{i+}$, then $\overline{x}\in
X_{q1}^{i}$ and $q_{0}(\overline{x})=\min_{x\in X_{i1}}q_{0}(x)=\max
_{\lambda\in Y_{1}^{i+}}D(\lambda)=D(\overline{\lambda})$. If $\overline
{\lambda}\in\mathbb{R}_{+}^{m}\cap Y^{-}$ then $q_{0}(\overline{x})=\max_{x\in
X_{i}}q_{0}(x)=\max_{\lambda\in\mathbb{R}_{+}^{m}\cap Y^{-}}D(\lambda
)=D(\overline{\lambda})$. }

\smallskip Clearly, the conclusion of \cite[Th.~4]{Gao:03a} follows from
Proposition \ref{p1i} (ii) and (i).

Let us look at \cite[Th.~6]{Gao:03a}. Because $(\overline{x},\overline
{\lambda})$ is a KKT point of $L$ with $\overline{\lambda}\in Y_{q1}%
^{i+}\subset Y^{i+}$, (\ref{r-minmaxqi}) holds. Moreover, because
$\overline{\lambda}_{1}>0$, it follows that $q_{1}(\overline{x})=0$, and so
$\overline{x}\in X_{q1}^{i}$ $(\subset X_{i})$, and so the first assertion of
\cite[Th.\ 6]{Gao:03a} holds, but (\ref{r-minmaxqi}) is stronger.

Consider now the particular case in which $b_{j}=0$ and $c_{j}=0$ for
$j\in\overline{2,m}$ (or, equivalently, $m=1$); in this case the preceding
problem becomes a \textquotedblleft quadratic programming problem over a
sphere\textquotedblright, considered in \cite[Sect.~6]{Gao:03a}. Assume that
$Y^{-}\ni\overline{\lambda}=\overline{\lambda}_{1}>0$. Then $\nabla
D(\overline{\lambda})=0$, and so $\overline{x}\in X_{e}$. Using Proposition
\ref{p1} we get
\[
\max_{x\in X_{i}}q_{0}(x)\geq q_{0}(\overline{x})=\max_{x\in X_{e}}%
q_{0}(x)=\min_{\lambda\in Y^{-}}D(\lambda)=\min_{\lambda\in
Y_{\operatorname{col}}^{-}}D(\lambda)=D(\overline{\lambda}),
\]
which does not agree with the second assertion of \cite[Th.~6]{Gao:03a}
because $\mathbb{R}_{+}\cap Y^{-}\subset Y^{-}\subset Y_{\operatorname{col}%
}^{-}$.

\begin{example}
\label{ex-gao03}Let $n=1$, $q_{0}(x)=-\tfrac{1}{2}(x^{2}+x)$ and
$q_{1}(x)=\tfrac{1}{2}(x^{2}-1)$. It follows that $X_{e}=\{-1,1\}$,
$X_{i}=[-1,1]$, $Y^{+}=(1,\infty)=Y^{i+}$ and $Y^{-}=(-\infty,1)\supset
\lbrack0,1)=\mathbb{R}_{+}\cap Y^{-}$. In this case we have that
$A(\lambda)=\lambda-1$, $b(\lambda)=\tfrac{1}{2}$, $c(\lambda)=-\frac{\lambda
}{2}$, $L(x,\lambda)=\frac{\lambda-1}{2}x^{2}-\tfrac{1}{2}x-\frac{\lambda}{2}%
$, $\nabla L(x,\lambda)=\left(  (\lambda-1)x-\tfrac{1}{2},\tfrac{1}{2}%
x^{2}-\tfrac{1}{2}\right)  $, $\nabla L(x,\lambda)=0$ $\Leftrightarrow$
$(x,\lambda)\in\left\{  (-1,\tfrac{1}{2}),(1,\tfrac{3}{2})\right\}  $,
$D(\lambda)=\frac{1}{8(1-\lambda)}-\frac{\lambda}{2}$. For $(\overline
{x},\overline{\lambda})=(1,\tfrac{3}{2})$ we have that
\[
q_{0}(\overline{x})=\min_{x\in X_{i}}q_{0}(x)=\max_{\lambda\in
Y_{\operatorname{col}}^{i+}}D(\lambda)=D(\overline{\lambda}),
\]
which confirms the second assertion of \cite[Th.~6]{Gao:03a}, while for
$(\overline{x},\overline{\lambda})=(-1,\tfrac{1}{2})$ we have that
\[
\tfrac{1}{8}=\max_{x\in\lbrack-1,1]}q_{0}(x)>0=q_{0}(-1)=\max_{x\in
\{-1,1\}}q_{0}(x)=\min_{\lambda\in\lbrack0,1)}D(\lambda)=D(\tfrac{1}{2}%
)<\sup_{\lambda\in\lbrack0,1)}D(\lambda)=\infty.
\]
This shows that the third assertion of \cite[Th.~6]{Gao:03a} is false.
\end{example}

Of course, in \cite[Th.~6]{Gao:03a} there is no need to assume $A$ (i.e.\ our
$A_{0}$) \textquotedblleft has at least one negative
eigenvalue\textquotedblright; probably this hypothesis was added in order
problem $(\mathcal{P}_{\lambda})$ be not a convex one.

\medskip

The problems considered by DY Gao in his survey papers \cite[Sect.~4]{Gao:08}
and \cite[Sect.~4]{Gao:09} (which are almost the same) refer to
\textquotedblleft box constrained problem\textquotedblright\ (\cite{Gao:07},
\cite{GaoRua:10}), \textquotedblleft integer programming\textquotedblright%
\ (\cite{FanGaoSheWu:08}, \cite{Gao:07}, \cite{GaoRua:10},
\cite{WanFanGaoXin:08}), \textquotedblleft mixed integer programming with
fixed charge\textquotedblright\ (\cite{GaoRuaShe:10}) and \textquotedblleft
quadratic constraints\textquotedblright\ (\cite{GaoRuaShe:09}). In these
survey papers the results are stated without proofs and their statements are
generally different from the corresponding ones in the papers mentioned above;
even more, for some results, the statements are different in the two survey
papers, even if the wording (text) is almost the same. We shall mention those
results from \cite[Sect.~4]{Gao:08} and/or \cite[Sect.~4]{Gao:09} which have
not equivalent statements in other papers.

It seems that the first paper dedicated completely to quadratic problems with
quadratic equality constraints using CDT is \cite{FanGaoSheWu:08}, even if
\cite{Gao:07} was published earlier; note that \cite{FanGaoSheWu:08} is cited
in \cite{Gao:07} as Ref.~6 with a slightly different title (see also Ref. Fang
SC, Gao DY, Sheu RL, Wu SY (2007a) in \cite{GaoRua:08}).

\medskip

The problems considered by Fang, Gao, Sheu and Wu in \cite{FanGaoSheWu:08} are
of type $(P_{e})$ with $m=n$. Setting $e_{j}:=(\delta_{jk})_{k\in
\overline{1,n}}\in\mathbb{R}^{n}$, one has $A_{j}:=2\operatorname*{diag}e_{j}%
$, $b_{j}:=e_{j}$, $c_{j}:=0$ for $j\in\overline{1,n}$. Of course,
$X_{e}=\{0,1\}^{n}$.

\smallskip Theorem 1 in \cite[Th.~1]{FanGaoSheWu:08} asserts: \emph{Let
$\overline{\lambda}\in Y_{0}\cap\mathbb{R}_{++}^{n}$ be a critical point of
$D$ and $\overline{x}:=x(\overline{\lambda})$. Then $\overline{x}$ is a KKT
point for problem $(P_{e})$ and $q_{0}(\overline{x})=D(\overline{\lambda})$. }

\smallskip Theorem 2 in \cite[Th.~1]{FanGaoSheWu:08} asserts: \emph{Let
$\overline{\lambda}\in Y_{0}\cap\mathbb{R}_{--}^{n}$ be a critical point of
$D$ and $\overline{x}:=x(\overline{\lambda})$. Then $\overline{x}$ is a KKT
point for the problem $(\mathcal{P}_{\max})$ of maximizing $q_{0}$ on $X_{e}$
and $q_{0}(\overline{x})=D(\overline{\lambda})$. }

\smallskip Theorem 3 in \cite[Th.~1]{FanGaoSheWu:08} asserts: \emph{Let
$\overline{\lambda}\in Y_{0}$ be a critical point of $D$ and $\overline
{x}:=x(\overline{\lambda})$. }

\emph{(a) If $\overline{\lambda}\in\mathcal{S}_{\natural}^{+}:=Y^{+}%
\cap\mathbb{R}_{++}^{n}$, then $q_{0}(\overline{x})=\min_{x\in X_{e}}%
q_{0}(x)=\max_{\lambda\in\mathcal{S}_{\natural}^{+}}D(\lambda)=D(\overline
{\lambda})$. }

\emph{(b) If $\overline{\lambda}\in\mathcal{S}_{\natural}^{-}:=Y^{-}%
\cap\mathbb{R}_{++}^{n}$, then in a neighborhood $\mathcal{X}_{0}%
\times\mathcal{S}_{0}\subset X_{e}\times\mathcal{S}_{\natural}^{-}$ of
$(\overline{x},\overline{\lambda})$, $q_{0}(\overline{x})=\min_{x\in
\mathcal{X}_{0}}q_{0}(x)=\min_{\lambda\in\mathcal{S}_{0}}D(\lambda
)=D(\overline{\lambda})$. }

\emph{(c) If $\overline{\lambda}\in\mathcal{S}_{\flat}^{-}:=Y^{-}%
\cap\mathbb{R}_{--}^{n}$, then $q_{0}(\overline{x})=\max_{x\in X_{e}}%
q_{0}(x)=\min_{\lambda\in\mathcal{S}_{\flat}^{-}}D(\lambda)=D(\overline
{\lambda})$. }

\emph{(d) If $\overline{\lambda}\in\mathcal{S}_{\flat}^{+}:=Y^{+}%
\cap\mathbb{R}_{--}^{n}$, then in a neighborhood $\mathcal{X}_{0}%
\times\mathcal{S}_{0}\subset X_{e}\times\mathcal{S}_{\flat}^{+}$ of
$(\overline{x},\overline{\lambda})$, $q_{0}(\overline{x})=\max_{x\in
\mathcal{X}_{0}}q_{0}(x)=\max_{\lambda\in\mathcal{S}_{0}}D(\lambda
)=D(\overline{\lambda})$. }

\medskip

Using Proposition 4 for $\overline{\lambda}\in Y_{0}$ with $\nabla
D(\overline{\lambda})=0$ we have: (i) $q_{0}(\overline{x})=D(\overline
{\lambda})$ without supplementary conditions on $\overline{\lambda};$ (ii)
because $\mathcal{S}_{\natural}^{+}\subset Y^{+}$, Eq.~(\ref{r-minmaxqe}) is
stronger than the minmax relation in (a); (iii) because $\mathcal{S}_{\flat
}^{-}\subset Y^{-}$, Eq.~(\ref{r-maxminqe}) is stronger than the maxmin
relation in (c); (iii) because $q_{0}$ is locally constant on $X_{e}$, (b) and
(d) are true but their conclusions are much weaker then those provided by
Eq.~(\ref{r-maxminqe}) and Eq.~(\ref{r-minmaxqe}), respectively.

\medskip

The quadratic problems $(\mathcal{P}_{b})$ considered by Gao in \cite[Th.~4]%
{Gao:07} is of type $(P_{e})$ in which $m\geq n$, $q_{0}(x):=-\tfrac{1}%
{2}\left\Vert Ax-c\right\Vert ^{2}$ for some $A\in\mathbb{R}^{p\times n}$ and
$c\in\mathbb{R}^{p}$, $A_{j}:=\operatorname*{diag}e_{j}$, $b_{j}:=0$,
$c_{j}:=-\tfrac{1}{2}$ for $j\in\overline{1,n}$, $A_{j}=0$ for $j\in
\overline{n+1,m}$; hence $X_{e}\subset\{-1,1\}^{n};$ problem $(\mathcal{P}%
_{bo})$ is $(\mathcal{P}_{b})$ in the case $m=n$. The problem of maximizing
$D$ on $\mathcal{S}_{b}:=Y_{0}\cap(\mathbb{R}_{++}^{n}\times\mathbb{R}^{m-n})$
is denoted by $(\mathcal{P}_{b}^{d})$ in the general case, and by
$(\mathcal{P}_{bo}^{d})$ for $m=n$ (when $\mathcal{S}_{b}:=Y_{0}\cap
\mathbb{R}_{++}^{n}$).

\smallskip Theorem 4 in \cite{Gao:07} asserts: \emph{Let $\overline{\lambda
}\in\mathcal{S}_{b}$ be \textquotedblleft a critical point of $(\mathcal{P}%
_{b}^{d})$\textquotedblright\ and $\overline{x}:=x(\overline{\lambda})$. Then
$\overline{x}$ \textquotedblleft is a critical point of $(\mathcal{P}_{b})$"
and $q_{0}(\overline{x})=D(\overline{\lambda})$. Moreover, if $\overline
{\lambda}\in\mathcal{S}_{b}^{+}:=Y^{+}\cap(\mathbb{R}_{++}^{n}\times
\mathbb{R}^{m-n})$, then $q_{0}(\overline{x})=\min_{x\in X_{e}}q_{0}%
(x)=\max_{\lambda\in\mathcal{S}_{b}^{+}}D(\lambda)=D(\overline{\lambda})$. }

\smallskip Corollary 2 in \cite{Gao:07} asserts: \emph{Let $\overline{\lambda
}\in\mathcal{S}_{b}$ be \textquotedblleft a KKT point the canonical dual
problem $(\mathcal{P}_{bo}^{d})$\textquotedblright\ and $\overline
{x}:=x(\overline{\lambda})$. Then $\overline{x}$ \textquotedblleft is a KKT
point of the Boolean least squares problem $(\mathcal{P}_{bo})$". If
$\overline{\lambda}\in\mathcal{S}_{b}^{+}$, then $q_{0}(\overline{x}%
)=\min_{x\in X_{e}}q_{0}(x)=\max_{\lambda\in\mathcal{S}_{b}^{+}}%
D(\lambda)=D(\overline{\lambda})$. }

Unfortunately, it is not defined what is meant by critical points of problems
$(\mathcal{P}_{b}^{d})$ and $(\mathcal{P}_{b})$, respectively. However,
because $\mathcal{S}_{b}$ and $\mathcal{S}_{b}^{+}$ are open sets, by
\textquotedblleft critical point of $(\mathcal{P}_{b}^{d})$\textquotedblright%
\ one must mean \textquotedblleft critical point of $D$\textquotedblright; in
this situation\ the conclusions of \cite[Th.~4]{Gao:07}, less $\overline{x}$
\textquotedblleft is a critical point of $(\mathcal{P}_{b})$", are true, but
are much weaker than those provided by Proposition \ref{p1}. Similarly, in
\cite[Cor.~2]{Gao:07}, $\overline{\lambda}\in\mathcal{S}_{b}$ is
\textquotedblleft a KKT point the canonical dual problem $(\mathcal{P}%
_{bo}^{d})$\textquotedblright\ is equivalent to $\overline{\lambda}$ is a
\textquotedblleft critical point of $D$\textquotedblright.

\medskip

The difference between problems $(\mathcal{P}_{b})$ considered by Wang, Fang,
Gao and Xing in \cite[p.~215]{WanFanGaoXin:08} and \cite[Th.~4]{Gao:07} is
that in the former $q_{0}$ is a general quadratic function (hence
$X_{e}\subset\{-1,1\}^{n}$).

\smallskip Theorem 2.2 in \cite{WanFanGaoXin:08} asserts: \emph{Let
$\overline{\lambda}\in\mathcal{S}_{b}:=Y_{0}\cap(\mathbb{R}_{+}^{n}%
\times\mathbb{R}^{m-n})$ be a critical point of $D$ and $\overline
{x}:=x(\overline{\lambda})$. Then $\overline{x}$ is a KKT point of $(P_{e})$
and $q_{0}(\overline{x})=D(\overline{\lambda})$. }

\smallskip Theorem 2.3 in \cite{WanFanGaoXin:08} asserts: \emph{Let
$\overline{\lambda}\in\mathcal{S}_{b}^{+}:=Y^{+}\cap(\mathbb{R}_{+}^{n}%
\times\mathbb{R}^{m-n})$ be a critical point of $D$ and $\overline
{x}:=x(\overline{\lambda})$. Then $q_{0}(\overline{x})=\min_{x\in X_{e}}%
q_{0}(x)=\max_{\lambda\in\mathcal{S}_{b}^{+}}D(\lambda)=D(\overline{\lambda}%
)$. }

\medskip

Theorems 3.2 and 3.3 from \cite{WanFanGaoXin:08} are the versions of Theorems
2.2 and 2.3 for $n=m$, respectively. Of course, the conclusions of Theorems
2.2 and 2.3 are valid replacing $\mathcal{S}_{b}$ and $\mathcal{S}_{b}^{+}$ by
$Y_{0}$ and $Y^{+}$, respectively.

\medskip

The general quadratic problem with inequality constraints $(P_{i})$ is
considered by Gao in \cite{Gao:08} and \cite{Gao:09}. In the sequel, the
Moore--Penrose generalized inverse of $F\in\mathfrak{M}_{n}$ is denoted by
$F^{\dag}$ or $F^{+}$, as in the corresponding cited papers authored by Gao
and his collaborators.

\smallskip Theorem 7 in \cite{Gao:08} and Theorem 10 in \cite{Gao:09} assert:
\emph{Let $\overline{\lambda}\in Y_{\operatorname{col}}^{i+}$ be a a solution
of problem $(\mathcal{P}_{q}^{d})$ of maximizing $D$ on $Y_{\operatorname{col}%
}^{i}$ and $\overline{x}:=\left[  A(\overline{\lambda})\right]  ^{\dag
}(\overline{\lambda})$. Then $\overline{x}$ is a KKT point of $(P_{i})$ and
$q_{0}(\overline{x})=D(\overline{\lambda})$. If $A(\overline{\lambda}%
)\succeq0$ then $\overline{\lambda}$ is a global maximizer of the problem
$(\mathcal{P}_{q}^{d})$ and $\overline{x}$ is a global minimizer of\ $(P_{i}%
)$. If $A(\overline{\lambda})\prec0$, then $\overline{x}$ is a local minimizer
(or maximizer) of\ $(P_{i})$ if and only if $\overline{\lambda}$ is a local
minimizer (or maximizer) of $D$ on $Y_{\operatorname{col}}^{i+}$.} \smallskip

\smallskip The \textquotedblleft box constrained problem" $(\mathcal{P}_{b})$
considered by Gao in \cite[Th.~3]{Gao:08} and \cite{Gao:09} is of type
$(P_{i})$ in which $m=n$, $A_{j}:=2\operatorname*{diag}e_{j}$, $b_{j}:=0$,
$c_{j}:=-1$ for $j\in\overline{1,n}$; hence $X_{i}=[-1,1]^{n}$.

\smallskip Theorem 3 in \cite{Gao:08} asserts: \emph{Let $\overline{\lambda
}\in Y_{\operatorname{col}}^{i+}$ be a critical point of $D$\ and
$\overline{x}:=\left[  A(\overline{\lambda})\right]  ^{\dag}b(\overline
{\lambda})$. Then $\overline{x}$ is a KKT point of $(P_{i})$ and
$q_{0}(\overline{x})=D(\overline{\lambda})$. Moreover, if $A(\overline
{\lambda})\succeq0$ then $q_{0}(\overline{x})=\min_{x\in X_{i}}q_{0}%
(x)=\max_{\lambda\in Y_{\operatorname{col}}^{i+}}D(\lambda)=D(\overline
{\lambda})$. If $A(\overline{\lambda})\prec0$, then on a neighborhood
$\mathcal{X}_{o}\times\mathcal{S}_{o}$ of $(\overline{x},\overline{\lambda})$
we have either $q_{0}(\overline{x})=\min_{x\in\mathcal{X}_{0}}q_{0}%
(x)=\min_{\lambda\in\mathcal{S}_{0}}D(\lambda)=D(\overline{\lambda})$, or
$q_{0}(\overline{x})=\max_{x\in\mathcal{X}_{0}}q_{0}(x)=\max_{\lambda
\in\mathcal{S}_{0}}D(\lambda)=D(\overline{\lambda})$.}

\smallskip The only difference between \cite[Th.~3]{Gao:08} and \cite[Th.~5]%
{Gao:09} is that in the latter the case $A(\overline{\lambda})\prec0$ is missing.

\smallskip

Probably, the intention was to take $\overline{\lambda}\in
Y_{\operatorname{col}}^{i}$ instead of $\overline{\lambda}\in
Y_{\operatorname{col}}^{i+}$ in the first assertions of \cite[Ths.~3,
7]{Gao:08} and \cite[Ths.~5, 10]{Gao:09}; in fact, there is not $\overline
{\lambda}\in Y_{\operatorname{col}}^{i+}$ such that $A(\overline{\lambda
})\prec0!$

It is not clear how the criticality of $D$ at $\lambda\in
Y_{\operatorname{col}}\setminus Y_{0}$ is defined in \cite[Th.~3]{Gao:08} and
\cite[Th.~5]{Gao:09}.

Let us assume that $\overline{\lambda}\in Y_{0}$ is a critical point of $D$ in
the mentioned results from \cite{Gao:08} and \cite{Gao:09}; in this situation
\cite[Th.~3]{Gao:08} is a particular case of \cite[Th.~7]{Gao:08}. Then
$\overline{x}$ is a KKT point of $(P_{i})$ iff $\overline{\lambda}%
\in\mathbb{R}_{+}^{n};$ assuming moreover that $A(\overline{\lambda})\succeq
0$, the conclusion of the second assertion of \cite[Th.~7]{Gao:08} is true.
However, in the case $A(\overline{\lambda})\prec0$ the conclusions of
\cite[Ths.~3, 7]{Gao:08} are false, as the next example shows.

\begin{example}
\label{ex-gao08} Consider $n:=m:=2$, $A_{0}:=\left[
\begin{array}
[c]{cc}%
-1 & 1\\
1 & -3
\end{array}
\right]  $, $A_{1}:=\operatorname*{diag}e_{1}$, $A_{2}:=\operatorname*{diag}%
e_{2}$, $b_{0}:=(0,-1)^{T},b_{1}:=b_{2}:=0$, $c_{1}:=c_{2}:=-\frac{1}{2}$.
Then $A(\lambda)=A_{0}+\lambda_{1}A_{1}+\lambda_{2}A_{2}$, $b(\lambda)=b_{0}$,
$c(\lambda)=-\tfrac{1}{2}(\lambda_{1}+\lambda_{2})$. We have that
$Y_{\operatorname{col}}=Y_{0}=\{(\lambda_{1},\lambda_{2})\in\mathbb{R}^{2}%
\mid(\lambda_{1}-1)(\lambda_{2}-3)\neq1\}$. The critical points $(\overline
{x},\overline{\lambda})$ of $L$ are: $\left(  (-1,-1)^{T},(0,3)^{T}\right)  $,
$\left(  (-1,1)^{T},(2,3)^{T}\right)  $, $\left(  (1,-1)^{T},(2,5)^{T}\right)
$, $\left(  (1,1)^{T},(0,1)^{T}\right)  $. Applying Proposition \ref{p1i} we
obtain that $\overline{x}:=(1,-1)^{T}$ is the global minimizer of $q_{0}$ on
$X_{i}=[0,1]^{2}$ and $\overline{\lambda}:=(2,5)^{T}$ is the global maximizer
of $D$ on $Y_{\operatorname{col}}^{i}=Y^{i}=\{(\lambda_{1},\lambda_{2}%
)\in\mathbb{R}^{2}\mid\lambda_{1}>2$, $(\lambda_{1}-2)(\lambda_{2}-3)>1\}$.

Take now $(\overline{x},\overline{\lambda}):=\left(  (1,1)^{T},(0,1)^{T}%
\right)  ;$ we have that $\overline{\lambda}\in\mathbb{R}_{+}^{2}$ and
$A(\overline{\lambda})\prec0$. From Proposition \ref{lem-pd}~(iv), we have
that $\overline{\lambda}$ is a global minimizer of $D$ on
$Y_{\operatorname{col}}^{-}=Y^{-}=\{(\lambda_{1},\lambda_{2})\in\mathbb{R}%
^{2}\mid\lambda_{1}<2$, $(\lambda_{1}-2)(\lambda_{2}-3)>1\}$. Assuming that
$\overline{\lambda}$ is a local maximizer of $D$, because $D$ is convex on
$Y_{\operatorname{col}}^{-}=Y^{-}$, $D$ is constant on an open neighborhood
$U\subset Y^{-}$ of $\overline{\lambda}$, and so $\nabla D(\lambda)=0$ for
$\lambda\in U;$ taking into account (\ref{r-cppdL}), this is a contradiction.
Observe that $\overline{x}=(1,1)$ is not a local minimizer of $q_{0}$ on
$X_{i}$. Indeed, take $x:=(1-u,1)\in X_{i}$ for $u\in(0,2);$ then
$q_{0}(x)=-\tfrac{1}{2}u^{2}<0=q_{0}(\overline{x})$, proving that
$\overline{x}$ is not a local minimum of $q_{0}$ on $X_{i}$.
\end{example}

Gao and Sherali in \cite[Th.\ 8.16]{GaoShe:09} (attributed to \cite{Gao:05})
assert: \emph{Suppose that $m=1$, $A_{1}>0$, $b_{1}=0$, $c_{1}<0$. Let
$\overline{\lambda}\in Y^{i}$ be a critical point of $D$ and $\overline
{x}:=x(\overline{\lambda})$. If $\overline{\lambda}\in Y^{i+}$, then
$\overline{x}$ is a global minimizer of $q_{0}$ on $X_{i}$. If $\overline
{\lambda}\in\mathbb{R}_{+}\cap Y^{-}$ then $\overline{\lambda}$ is a local
minimizer of $q_{0}$ on $X_{i}$. }

\smallskip As in the case of \cite[Th.~6]{Gao:03a} above, the first assertion
of \cite[Th.~8.16]{GaoShe:09} follows from Proposition \ref{p1i}. However, the
second assertion of \cite[Th.~8.16]{GaoShe:09} is false as the next example shows.

\begin{example}
\emph{(see \cite[Ex.~1]{VoiZal:10})}\label{ex-th16-gs} Consider $n:=2$,
$m:=1$, $A_{0}:=\left[
\begin{array}
[c]{cc}%
-2 & -1\\
-1 & -3
\end{array}
\right]  $, $A_{1}:=I_{2}$, $b_{0}:=(-1,-1)^{T}$, $b_{1}:=0$, $c_{1}%
:=-\frac{1}{2}$. Then $D(\lambda)=-\frac{1}{2}\lambda-\frac{1}{2}%
\frac{2\lambda-3}{\lambda^{2}-5\lambda+5}$ and $D^{\prime}(\lambda)=-\frac
{1}{2}\frac{\left(  \lambda-2\right)  ^{2}}{\left(  \lambda^{2}-5\lambda
+5\right)  ^{2}}\left(  \lambda-1)(\lambda-5\right)  $. Hence the set of
critical points of $D$ is $\{1,2,5\}\subset\mathbb{R}_{+}$. For $\overline
{\lambda}=1$ we have that $A(\overline{\lambda})=\left(
\begin{array}
[c]{cc}%
-1 & -1\\
-1 & -2
\end{array}
\right)  \prec0$ and $\overline{x}=x(\overline{\lambda})=(1,0)^{T}$. Since
$X_{i}=\{(\cos t,\sin t)^{T}\mid t\in(-\pi,\pi]\}$ and
\[
q_{0}((\cos t,\sin t)^{T})=-(3+\cos t-2\sin t)\sin^{2}\tfrac{1}{2}t\leq
(\sqrt{5}-3)\sin^{2}\tfrac{1}{2}t<0=q_{0}(\overline{x})
\]
for all $t\in(-\pi,\pi]\setminus\{0\}$, we have that $\overline{x}$ is the
unique global maximizer of $q_{0}$ on $X_{i}$, in contradiction with the
second assertion of \cite[Th.~8.16]{GaoShe:09}.
\end{example}


The problem considered by Zhang, Zhu and Gao in \cite{ZhaZhuGao:09} is of type
$(P_{i})$ in which $m\geq n$, $A_{j}:=\operatorname*{diag}e_{j}$, $b_{j}:=0$,
$c_{j}\leq0$ for $j\in\overline{1,n}$, $A_{j}=0$ for $j\in\overline{n+1,m}$.

\smallskip Theorem 1 in \cite{ZhaZhuGao:09} asserts: \emph{Let $\overline
{\lambda}\in Y^{i}$ be a KKT point of $D$ and $\overline{x}:=x(\overline
{\lambda})$. Then $\overline{x}$ is a KKT point of $(P_{i})$ and
$q_{0}(\overline{x})=D(\overline{\lambda})$. }

\smallskip Theorem 2 in \cite{ZhaZhuGao:09} asserts: \emph{Let $\overline
{\lambda}\in Y^{i}$ be a KKT point of $D$ and $\overline{x}:=x(\overline
{\lambda})$. If $\overline{\lambda}\in Y^{i+}$, then $\overline{\lambda}$
\textquotedblleft is a global maximizer of\textquotedblright\ $D$ on $Y^{i+}$
\textquotedblleft if and only if the vector" $\overline{x}$ \textquotedblleft
is a global minimizer of\textquotedblright\ $(P_{i})$ on $X_{i}$, and
$q_{0}(\overline{x})=\min_{x\in X_{i}}q_{0}(x)=\max{}_{\lambda\in Y^{i+}%
}D(\lambda)=D(\overline{\lambda})$. If $\overline{\lambda}\in\mathbb{R}%
_{+}^{m}\cap Y^{-}$, \textquotedblleft then in a neighborhood $\mathcal{X}%
_{0}\times S_{0}\subset$\textquotedblright$X_{i}\times(\mathbb{R}_{+}^{m}\cap
Y^{-})$ of $(\overline{x},\overline{\lambda})$, \textquotedblleft we have that
either\textquotedblright\ $q_{0}(\overline{x})=\min_{x\in\mathcal{X}_{0}}%
q_{0}(x)=\min_{\lambda\in S_{0}}D(\lambda)=D(\overline{\lambda})$, or
$q_{0}(\overline{x})=\max_{x\in\mathcal{X}_{0}}q_{0}(x)=\max{}_{\lambda\in
S_{0}}D(\lambda)=D(\overline{\lambda})$. }

\medskip

Clearly, \cite[Th.~1]{ZhaZhuGao:09} and the conclusion of \cite[Th.~2]%
{ZhaZhuGao:09} in the case $\overline{\lambda}\in Y^{i+}$ follow from
Proposition \ref{p1i}. As shown in \cite[Ex.\ 2]{Zal:11} and Example
\ref{ex-gao08}, each of the alternative conclusions of \cite[Th.~2]%
{ZhaZhuGao:09} in the case $\overline{\lambda}\in\mathbb{R}_{+}^{m}\cap Y^{-}$
is false. Observe that \cite{GaoRuaShe:09} is cited in \cite{ZhaZhuGao:09} as
a paper to appear, but not in connection with the previous result.

\medskip

The problem $(P_{i})$ is considered also by Gao, Ruan and Sherali in
\cite[p.~486]{GaoRuaShe:09}; the problem of maximizing $D$ on
$Y_{\operatorname{col}}^{i}$ is denoted by $(\mathcal{P}_{q}^{d})$.

\smallskip Theorem 4 in \cite{GaoRuaShe:09} asserts: \emph{Let $\overline
{\lambda}\in Y_{\operatorname{col}}^{i}$ be a critical point of $(\mathcal{P}%
_{q}^{d})$ and $\overline{x}:=\left[  A(\overline{\lambda})\right]
^{+}b(\overline{\lambda})$. Then $\overline{x}$ is a KKT point of $(P_{i})$
and $q_{0}(\overline{x})=D(\overline{\lambda})$. If $\overline{\lambda}\in
Y_{\operatorname{col}}^{i+}$, then $q_{0}(\overline{x})=\min_{x\in X_{i}}%
q_{0}(x)=\max_{\lambda\in Y_{\operatorname{col}}^{i+}}D(\lambda)=D(\overline
{\lambda})$. If $\overline{\lambda}\in Y^{i+}$ then $\overline{\lambda}$
\textquotedblleft is a unique global maximizer of $(\mathcal{P}_{q}^{d})$ and
the vector $\overline{x}$ is a unique global minimizer of $(P_{i}%
)$\textquotedblright. If $\overline{\lambda}\in\mathbb{R}_{+}^{m}\cap Y^{-}$,
then $\overline{\lambda}$ \textquotedblleft is a local minimizer
of\textquotedblright\ $D$ \textquotedblleft on the neighborhood $S_{o}\subset
$\textquotedblright$\mathbb{R}_{+}^{m}\cap Y^{-}$ \textquotedblleft if and
only if $\overline{x}$ is a local minimizer of\textquotedblright\ $q_{0}$
\textquotedblleft on the neighborhood $X_{o}\subset$\textquotedblright$X_{i}$,
i.e., $q_{0}(\overline{x})=\min_{x\in\mathcal{X}_{o}}q_{0}(x)=\min_{\lambda
\in\mathcal{S}_{o}}D(\lambda)=D(\overline{\lambda})$.}

\smallskip As noticed before Lemma \ref{lem-im}, $Y_{\operatorname{col}}$ is
not open in general, so it is not possible to speak about the
differentiability of $D$ at $\lambda\in Y_{\operatorname{col}}\setminus Y_{0}%
$. As in \cite[Th.~4]{Gao:07}, it is not explained what is meant by critical
point of $(\mathcal{P}_{q}^{d})$; we interpret it as being a critical point of
$D$. With the above interpretation for \textquotedblleft critical point of
$(\mathcal{P}_{q}^{d})$\textquotedblright, we agree with the first two
assertions of \cite[Th.~4]{GaoRuaShe:09}. However, the third assertion of
\cite[Th.~4]{GaoRuaShe:09} that $\overline{\lambda}$ is the unique global
maximizer of $(\mathcal{P}_{q}^{d})$ provided that $\overline{\lambda}\in
Y^{i+}$ is false, as seen in Example \ref{ex-qi-grs} below. The same example
shows that the fourth assertion of \cite[Th.~4]{GaoRuaShe:09} is false, too;
another counterexample is provided by Example \ref{ex-gao08}.

\begin{example}
\label{ex-qi-grs}Let us take $n=m=2$, $q_{0}(x,y):=xy-x$, and $q_{1}%
(x,y):=-q_{2}(x,y):=\tfrac{1}{2}\left(  x^{2}+y^{2}-1\right)  $ for
$(x,y)\in\mathbb{R}^{2}$. Clearly, the problems $(P_{e})$ for $(q_{0},q_{1})$
and $(P_{i})$ for $(q_{0},q_{1},q_{2})$ are equivalent in the sense that they
have the same objective functions and the same feasible sets (hence the same
solutions). Denoting by $L^{e}$, $A^{e}$, $b^{e}$, $c^{e}$, $D^{e}$ and
$L^{i}$, $A^{i}$, $b^{i}$, $c^{i}$, $D^{i}$ the functions associated to
problems $(P_{e})$ and $(P_{i})$ mentioned above, we get: $L^{e}%
(x,y,\lambda)=xy-x+\tfrac{\lambda}{2}\left(  x^{2}+y^{2}-1\right)  $,
$A^{e}(\lambda)=\left(
\begin{array}
[c]{ll}%
\lambda & 1\\
1 & \lambda
\end{array}
\right)  $, $b^{e}(\lambda)=(1,0)^{T}$, $c^{e}(\lambda)=-\tfrac{1}{2}\lambda$,
$Y_{\operatorname{col}}=Y_{0}=\mathbb{R}\setminus\{-1,1\}$,
$Y_{\operatorname{col}}^{+}=-Y_{\operatorname{col}}^{-}=Y^{+}=-Y^{-}%
=(1,\infty)$, $D^{e}(\lambda)=\frac{-\lambda}{\lambda^{2}-1}-\tfrac{1}%
{2}\lambda$ [for the problem $(P_{e})$] and $L^{i}(x,y,\lambda_{1},\lambda
_{2})=L^{e}(x,y,\lambda_{1}-\lambda_{2})$, $A^{i}(\lambda_{1},\lambda
_{2})=A^{e}(\lambda_{1}-\lambda_{2})$, $b^{i}(\lambda_{1},\lambda_{2}%
)=b^{e}(\lambda_{1}-\lambda_{2})$, $c^{i}(\lambda_{1},\lambda_{2}%
)=c^{e}(\lambda_{1}-\lambda_{2})$, $Y_{\operatorname{col}}^{i}=Y^{i}%
=\{(\lambda_{1},\lambda_{2})\in\mathbb{R}_{+}^{2}\mid\lambda_{1}-\lambda
_{2}\neq\pm1\}$, $Y_{\operatorname{col}}^{i+}=Y^{i+}=\{(\lambda_{1}%
,\lambda_{2})\in\mathbb{R}_{+}^{2}\mid\lambda_{1}-\lambda_{2}>1\}$,
$\mathbb{R}_{+}^{2}\cap Y^{-}=\{(\lambda_{1},\lambda_{2})\in\mathbb{R}_{+}%
^{2}\mid\lambda_{1}-\lambda_{2}<-1\}$, $D^{i}(\lambda_{1},\lambda_{2}%
)=D^{e}(\lambda_{1}-\lambda_{2})$.

The critical points of $L^{e}$ are $(0,1,0)$ and $\left(  \pm\sqrt
{3}/2,-1/2,\pm\sqrt{3}\right)  $. Using Proposition \ref{p1}, it follows that
$(\sqrt{3}/2,-1/2)$ is the unique global minimizer of $q_{0}$ on $X_{e}$ and
$\sqrt{3}$ is a global maximizer of $D^{e}$ on $Y_{\operatorname{col}}^{+}$
$(=Y^{+})$, while $(-\sqrt{3}/2,-1/2)$ is the unique global maximizer of
$q_{0}$ on $X_{e}$ and $-\sqrt{3}$ is a global minimizer of $D^{e}$ on
$Y_{\operatorname{col}}^{-}$ $(=Y^{-})$.

Note that $(\overline{x},\overline{y},\overline{\lambda}_{1},\overline
{\lambda}_{2})$ is a KKT point of $L^{i}$ iff $(\overline{x},\overline
{y},\overline{\lambda}_{1},\overline{\lambda}_{2})$ is a critical point of
$L^{i}$ with $(\overline{\lambda}_{1},\overline{\lambda}_{2})\in\mathbb{R}%
_{+}^{2}$, iff $(\overline{x},\overline{y},\overline{\lambda}_{1}%
-\overline{\lambda}_{2})$ is a critical point of $L^{e}$ with $(\overline
{\lambda}_{1},\overline{\lambda}_{2})\in\mathbb{R}_{+}^{2}$. Using Proposition
\ref{p1i}~(ii) we obtain that $(\sqrt{3}/2,-1/2)$ is the unique global
minimizer of $q_{0}$ on $X_{i}$ and any $(\overline{\lambda}_{1}%
,\overline{\lambda}_{2})\in\mathbb{R}_{+}^{2}$ with $\overline{\lambda}%
_{1}-\overline{\lambda}_{2}=\sqrt{3}$ is a global maximizer of $D^{i}$ on
$Y^{i+}$ $(=Y_{\operatorname{col}}^{i+})$, the latter assertion contradicting
the third assertion of \cite[Th.~4]{GaoRuaShe:09}. On the other hand, as seen
above, $(-\sqrt{3}/2,-1/2)$ is the unique global maximizer of $q_{0}$ on
$X_{e}=X_{i}$ and $(\sqrt{3},2\sqrt{3})\in\mathbb{R}_{+}^{2}\cap Y^{-}$ is a
global minimizer of $(\mathcal{P}_{q}^{d})$, contradicting the fourth
assertion of \cite[Th.~4]{GaoRuaShe:09}.
\end{example}

The problem considered by Lu, Wang, Xin and Fang in \cite{LuWanXinFan:10} is
of type $(P_{e})$ with $m=n$. More precisely, $A_{j}=2\operatorname*{diag}%
e_{j}$, $b_{j}:=e_{j}$, $c_{j}:=0$ for $j\in\overline{1,n};$ hence
$X_{e}=\{0,1\}^{n}$. One must emphasize the fact that the authors use the
usual Lagrangian, even if CDT is invoked.

\smallskip Theorem 2.2 (resp.\ Theorem 2.3) of \cite{LuWanXinFan:10} asserts:
\emph{If $\overline{\lambda}\in Y_{0}$ (resp.\ $\overline{\lambda}\in Y^{+}$)
is such that $\nabla D(\overline{\lambda})=0$ and $\overline{x}:=x(\overline
{\lambda})$, then $q_{0}(\overline{x})=D(\overline{\lambda})$ (resp.\ $q_{0}%
(\overline{x})=\min_{x\in X_{e}}q_{0}(x)$). }

\medskip

Gao and Ruan \cite{GaoRua:10} considered problems $(P_{e})$ and $(P_{i})$ when
$m=n$ and $A_{j}:=\operatorname*{diag}e_{j}$, $b_{j}:=0$, $c_{j}:=-\tfrac
{1}{2}$ for $j\in\overline{1,n}$. Of course, $X_{e}=\{-1,1\}^{n}$ and
$X_{i}=[-1,1]^{n}$. The problem of maximizing $D$ on $Y^{i+}$ is denoted by
$(\mathcal{P}^{d})$.

\smallskip Theorem 1 in \cite{GaoRua:10} (attributed to \cite{Gao:07})
asserts: \emph{\textquotedblleft If $\overline{\sigma}$ is a critical point
of\textquotedblright\ $D$, \textquotedblleft the vector $\overline{x}%
$"$:=x(\overline{\sigma})$ \textquotedblleft is a KKT point of" $(P_{i})$ and
$q_{0}(\overline{x})=D(\overline{\sigma})$. \textquotedblleft If the critical
point $\overline{\sigma}>0$, then the vector $\overline{x}$"$\in X_{e}$
\textquotedblleft is a local optimal solution of the integer programming
problem" $(P_{e})$. If $\overline{\sigma}\in Y^{i+}$, then $q_{0}(\overline
{x})=\min_{x\in X_{i}}q_{0}(x)=\max_{\sigma\in Y^{i+}}D(\sigma)=D(\overline
{\sigma})$. \textquotedblleft If the critical point $\overline{\sigma}\in
$"$Y^{i+}$ \textquotedblleft and $\overline{\sigma}>0$, then the vector
$\overline{x}$"$\in X_{e}$ \textquotedblleft is a global minimizer to the
integer programming problem" $(P_{e})$. If $\overline{\sigma}\in\mathbb{R}%
_{+}^{n}\cap Y^{-}$, \textquotedblleft then $\overline{\sigma}$ is a local
minimizer of $(\mathcal{P}^{d})$, the vector $\overline{x}$ is a local
minimizer of" $(P_{i})$, \textquotedblleft and on the neighborhood
$\mathcal{X}_{o}\times\mathcal{S}_{o}$ of $(\overline{x},\overline{\sigma})$,
$q_{0}(\overline{x})=\min_{x\in\mathcal{X}_{o}}q_{0}(x)=\min_{\sigma
\in\mathcal{S}_{o}}D(\sigma)=D(\overline{\sigma})$. }

\medskip

Concerning \cite[Th.~1]{GaoRua:10} we observe the following: In the first
assertion it is not clear if $\overline{\sigma}$ belongs to $\mathbb{R}%
_{+}^{n}$ or not; of course, $\overline{x}$ is not a KKT point of $(P_{i})$ if
$\overline{\sigma}\notin\mathbb{R}_{+}^{n}$. The second assertion is true
because $X_{e}$ is finite (without any condition on $\overline{\sigma}$). The
third assertion is false without assuming that $\overline{\sigma}$ is at least
a KKT point of $D$. The fourth assertion is true without assuming
$\overline{\sigma}>0$. The fifth assertion is false if $\overline{\sigma}>0$
and $\nabla D(\overline{\sigma})\neq0$.

\smallskip The main difference between \cite[Th.~1]{GaoRua:10} and the
conjunction of \cite[Th.~2 \& Th.~3]{GaoRua:10} is that in the latter $Y_{0}$
is replaced by $Y_{\operatorname{col}}$, but their statements are not more
clear. This is the reason for not analyzing them here.

\medskip

The problem considered by Gao, Ruan and Sherali in \cite{GaoRuaShe:10} is of
type $(P_{J})$ with $n=m=2k$ $(k\in\mathbb{N}^{\ast})$ and $J:=\overline
{k+1,n}$. In \cite{GaoRuaShe:10} $A_{0}$ is such that $(A_{0})_{ij}=0$ if
$\max\{i,j\}>k$, $A_{j}:=2\operatorname*{diag}e_{j}$ and $c_{j}:=0$ for
$j\in\overline{1,m}$, $b_{j}:=e_{j+k}$ for $j\in J^{c}$ $(=\overline{1,k})$
and $b_{j}:=e_{j}$ for $j\in J;$ moreover, $\mathcal{S}_{\natural
}:=Y_{\operatorname{col}}\cap(\mathbb{R}_{+}^{k}\times\mathbb{R}_{++}^{k})$
$(\subset Y_{\operatorname{col}}^{i}\subset Y_{\operatorname{col}}^{J})$,
$\mathcal{S}_{\natural}^{+}:=Y^{+}\cap\mathcal{S}_{\natural}$ $(\subset
Y^{i+}\subset Y^{J+})$, $\mathcal{S}_{\flat}:=Y_{\operatorname{col}}%
\cap(\mathbb{R}_{-}^{k}\times\mathbb{R}_{--}^{k})$, $\mathcal{S}_{\flat}%
^{-}:=Y^{-}\cap\mathcal{S}_{\flat}$ $(\subset Y^{i-}\subset Y^{J-})$.

\smallskip Theorem 1 of \cite{GaoRuaShe:10} asserts: \emph{Let $\overline
{\lambda}\in\mathcal{S}_{\natural}$ be a KKT point of $D$ and $\overline
{x}:=x(\overline{\lambda})$. Then $\overline{x}$ is feasible to the primal
problem $(P_{J})$ and $q_{0}(\overline{x})=L(\overline{x},\overline{\lambda
})=D(\overline{\lambda})$.}

\smallskip Theorem 2 in \cite{GaoRuaShe:10} asserts: \emph{Let $\overline
{\lambda}\in\mathcal{S}_{\natural}^{+}\cup\mathcal{S}_{\flat}^{-}$ be a
critical point of $D$ and $\overline{x}:=x(\overline{\lambda})$. If
$\overline{\lambda}\in\mathcal{S}_{\natural}^{+}$ then $q_{0}(\overline
{x})=\min_{x\in X_{J}}q_{0}(x)=\max_{\lambda\in\mathcal{S}_{\natural}^{+}%
}D(\lambda)=D(\overline{\lambda})$. If $\overline{\lambda}\in\mathcal{S}%
_{\flat}^{-}$ then $q_{0}(\overline{x})=\max_{x\in X_{J}}q_{0}(x)=\min
_{\lambda\in\mathcal{S}_{\flat}^{-}}D(\lambda)=D(\overline{\lambda})$. }

\smallskip In \cite[Th.~1]{GaoRuaShe:10} it is not clear what is meant by KKT
point of $D$ because $D$ is not differentiable for $\overline{\lambda}%
\in\mathcal{S}_{\natural}\setminus Y_{0}$. Propositions \ref{p1ei} and
\ref{p1eimax} confirm \cite[Th.~2]{GaoRuaShe:10}, but the conclusions of the
latter are much weaker than those of the former.

\medskip

The quadratic problems $(\mathcal{P}_{b})$ and $(\mathcal{P}_{bo})$ considered
by Ruan and Gao in \cite{RuaGao:17} (and \cite{RuaGao:16}) are those from
\cite{Gao:07}. The statement of \cite[Th.\ 5]{RuaGao:17} is that of
\cite[Cor.~2]{Gao:07} in which $\mathcal{S}_{b}$ is now $Y_{0}\cap\{\lambda
\in\mathbb{R}^{m}\mid\lambda_{j}\neq0$ $\forall j\in\overline{1,n}\}$,
$\mathcal{S}_{b}^{+}$ being the same, that is $Y^{+}\cap\mathbb{R}_{++}^{m}$.
The statement of \cite[Th.\ 6]{RuaGao:17} is that of \cite[Th.~4]{Gao:07} in
which \textquotedblleft a critical point of $(\mathcal{P}_{b}^{d}%
)$\textquotedblright\ is replaced by \textquotedblleft a KKT point of
$(\mathcal{P}_{b}^{d})$\textquotedblright.

\medskip

The quadratic problem considered by Ruan and Gao in \cite{RuaGao:17b} is of
type $(P_{J})$ in which $m>n$, and $\overline{1,n+1}\subset J$. In
\cite{RuaGao:17b} $A_{j}:=2\operatorname*{diag}e_{j}$, $b_{j}:=e_{j}$,
$c_{j}:=0$ for $j\in\overline{1,n}$, $A_{j}:=0$ for $j\in\overline{n+1,m};$
hence $X_{J}\subset\{0,1\}^{n}$. One considers $\mathcal{S}_{a}:=\{\lambda\in
Y^{J}\mid\lambda_{j}\neq0$ $\forall j\in J\}$ and $\mathcal{S}_{a}%
^{+}:=\{\lambda\in Y^{J+}\mid\lambda_{j}>0$ $\forall j\in J\}$.

\smallskip Theorem 3 of \cite{RuaGao:17b} asserts: \emph{Let $\overline
{\lambda}\in\mathcal{S}_{a}$ be a $J$-LKKT point of $D$ and $\overline
{x}:=x(\overline{\lambda})$. Then $\overline{x}$ is a $J$-LKKT point of
$(P_{J})$ and $q_{0}(\overline{x})=D(\overline{\lambda})$. }

\smallskip Theorem 4 in \cite{RuaGao:17b} asserts: \emph{Let $\overline
{\lambda}\in\mathcal{S}_{a}^{+}$ be a $J$-LKKT point of $D$ and $\overline
{x}:=x(\overline{\lambda})$. Then $q_{0}(\overline{x})=\min_{x\in X_{J}}%
q_{0}(x)=\max_{\lambda\in\mathcal{S}_{a}^{+}}D(\lambda)=D(\overline{\lambda}%
)$. }

\smallskip Clearly, \cite[Th.~3]{RuaGao:17b} is an immediate consequence of
Lemma \ref{lem-qperfdual}, while \cite[Th.~4]{RuaGao:17b} is a very particular
case of Proposition \ref{p1ei}.

\medskip

The quadratic problem considered by Gao in \cite{Gao:17}, \cite{Gao:18} and
\cite{GaoAli:18} is of type $(P_{J})$ in which $m=n+1$ and $J:=\overline{1,n}%
$. In these papers $A_{0}:=0$, $A_{j}:=2\operatorname*{diag}e_{j}$,
$b_{j}:=e_{j}$, $c_{j}:=0$ for $j\in J$, and $A_{n+1}:=0$; hence $X_{J}%
\subset\{0,1\}^{n}$.

\smallskip Theorem 2 of \cite{Gao:17} asserts: \emph{Let $\overline{\lambda
}\in Y^{J+}$ be a global maximizer of $D$ on $Y^{J+}$. Then $\overline
{x}:=x(\overline{\lambda})\in X_{J}$ and $q_{0}(\overline{x})=\min_{x\in
X_{J}}q_{0}(x)=\max_{\lambda\in Y^{J+}}D(\lambda)=D(\overline{\lambda})$. }

\smallskip The differences between \cite[Th.~2]{Gao:17} and \cite[Th.~2]%
{Gao:18} are: in the latter $b_{n+1}:=-c(u)\in\mathbb{R}_{-}^{n}$,
$c_{n+1}:=-V_{c}<0$, and $Y^{J+}$ is replaced by $\{\lambda\in Y^{J+}%
\mid\lambda_{n+1}>0\}$. The differences between \cite[Th.~2]{Gao:17} and
\cite[Th.~1]{GaoAli:18} are: in the latter $c_{n+1}:=-V_{c}<0$, and
$\min_{\rho\in\mathcal{Z}_{a}}P_{u}(\rho)$ is replaced by $\min_{\rho
\in\mathbb{R}^{n}}P_{u}(\rho);$ of course, $\min_{\rho\in\mathbb{R}^{n}}%
P_{u}(\rho)=-\infty$ if $c_{u}\neq0$. In all 3 papers there are provided
proofs of the mentioned results.

Using Proposition \ref{p1ei} (iii) in then context of \cite[Th.~2]{Gao:17} we
have that $\overline{\lambda}$ is a $J$-LKKT point of $D;$ using Proposition
\ref{p1ei} (ii) and (i) we get the conclusion of \cite[Th.~2]{Gao:17}.

\medskip

Yuan \cite{Yua:17} (the same as \cite{Yua:16}) considers problem $(P_{i})$ in
its general form.

\smallskip In \cite[p.~340]{Yua:17} one asserts: \emph{\textquotedblleft One
hard restriction is given" by $b_{0}\neq0$. \textquotedblleft The restriction
is very important to guarantee the uniqueness of a globally optimal solution
of" $(P_{i})$.}

\smallskip Theorem 1 of \cite{Yua:17} asserts: \emph{Let $\mathcal{Y}%
:=\{\sigma\in Y^{i}\mid x(\sigma)\in X_{i}\}\neq\emptyset$, and let
$(\mathcal{P}^{d})$ be the problem of maximizing $D$ on $\mathcal{Y}$. If
$\overline{\sigma}$ is a solution of $(\mathcal{P}^{d})$, then $\overline
{x}:=x(\overline{\sigma})$ is a solution of $(P_{i})$ and $q_{0}(\overline
{x})=D(\overline{\sigma})$. }

\smallskip Theorem 2 of \cite{Yua:17} asserts: \emph{Assume that ($C_{1}$)
$\sum_{k=0}^{m}A_{k}\succ0$, and ($C_{2}$) there exists $k\in\overline{1,m}$
such that $A_{k}\succ0$, $A_{0}+A_{k}\succ0$, and $\left\Vert D_{k}A_{0}%
^{-1}b_{0}\right\Vert >\left\Vert b_{k}^{T}D_{k}^{-1}\right\Vert
+\sqrt{\left\Vert b_{k}^{T}D_{k}^{-1}\right\Vert ^{2}+2|c_{k}|}$, where
$A_{k}=D_{k}^{T}D_{k}$ and $\left\Vert ^{\ast}\right\Vert $ is some vector
norm. Then problem $(\mathcal{P}^{d})$ has a unique non-zero solution
$\overline{\sigma}$ in the space $Y^{i+}$. }

\smallskip Counterexamples to both theorems of \cite{Yua:17} as well as for
the assertion on the \textquotedblleft hard restriction" $b_{0}\neq0$ from
\cite[p.~340]{Yua:17} are provided in \cite{Zal:18}.

\section{Conclusions}

\null
 -- We made a complete study of quadratic minimization problems
with quadratic equality and/or inequality constraints using the
method suggested by the canonical duality theory (CDT) introduced by
DY Gao. This method is based on the introduction of a dual function.
Our study uses only the usual Lagrangian associated to minimization
problems with equality and/or inequality constraints, without any
reference to CDT; CDT is presented (or, at least, referred) in all
the papers cited in Section \ref{sec4}.

-- As observed in Remark \ref{rem-adv}, it is more advantageous to use the
assertions (i) of Propositions \ref{p1}, \ref{p1ei}, \ref{p1i}, than the
second part of (ii) with $\overline{\lambda}\in Y_{0}$ because in versions (i)
one must know only the Lagrangian (hence only the data of the problems), and
this provides both $\overline{x}$ and $\overline{\lambda}$. Using $D$ could be
useful, possibly, if the number of constraints is much smaller than $n$.

-- As seen in Section \ref{sec4}, many results obtained by DY Gao and his
collaborators on quadratic optimization problems are not stated clearly, and
some of them are even false; some statements were made more clear in
subsequent papers, but we didn't observe some warning about the false
assertions. For the great majority of the correct assertions the use of the
usual direct method provides stronger versions.

-- Asking the strict positivity of the multipliers corresponding to nonlinear
constraints (but not only, as in \cite{RuaGao:17b}), is very demanding, even
for inequality constraints. Just observe that for $k$ equality constraints one
has $2^{k}$ distinct possibilities to get the feasible set, but at most one
could produce strictly positive multipliers.

\medskip

\textbf{Acknowledgement} We thank prof. Marius Durea for reading a previous
version of the paper and for his useful remarks.

\end{document}